\documentclass[12pt]{article}
\usepackage[english]{babel}
\usepackage{graphics}
\usepackage{amsfonts}
\usepackage{amssymb}
\usepackage{mathrsfs}
\usepackage{enumerate}
\usepackage{hyperref}
\usepackage{graphicx}

\newtheorem{corollary}{Corollary}[section]
\newtheorem{proposition}{Proposition}[section]
\newtheorem{example}{Example}[section]
\newtheorem{theorem}{Theorem}[section]

\newtheorem{remark}{Remark}[section]
\newtheorem{definition}{Definition}[section]
\def\QED{\hfill $\Box$\par\smallskip\noindent}

\def\scatola{\lower5pt\hbox{\vbox{\hrule\hbox{\vrule\kern2pt\vbox%
{\kern5pt\hbox{\mathsurround=0pt
}\kern2pt}\kern4pt\vrule}\hrule}}\ } 
\def\l{\lambda}

\def\convex{\setlength{{\baselineskip}{8pt}{\scriptsize{$\langle
\xi_1(\cdot),\cdot\rangle$}\normalsize}}}

\def\H{\mathbb{H}}

\def\R{\mathbb{R}}
\def\E{\mathbb{E}}
\def\convex{\scriptsize{$\langle \xi_1(\cdot),\cdot\rangle$}\, \normalsize}

\def\l{\lambda}
\def\g{\gamma}
\def\X{\mathbb{X}}

\def\0{\overline 0}

\def\psn{\par\smallskip\noindent}

\def\Cca{{\cal C}}

\def\Hca{{\cal H}}

\def\Pca{{\cal P}}

\def\Sca{{\cal S}}
\def\Tca{{\cal T}}

\def\G{\mathbf{G}}
\font\tenmsb=msbm10 \font\sevenmsb=msbm7 \font\fivemsb=msbm5
\newfam\msbfam
\textfont\msbfam=\tenmsb \scriptfont\msbfam=\sevenmsb
\scriptscriptfont\msbfam=\fivemsb
\font\teneufm=eufm10 \font\seveneufm=eufm7 \font\fiveeufm=eufm5
\newfam\eufmfam
\textfont\eufmfam=\teneufm \scriptfont\eufmfam=\seveneufm
\scriptscriptfont\eufmfam=\fiveeufm
\def\frak#1{{\fam\eufmfam\relax#1}}
\begin{document}
\author{A. Calogero\thanks{Dipartimento di Statistica,
Universit\`a degli Studi di Milano Bicocca, Via Bicocca degli
Arcimboldi 8, I-20126 Milano \tt(andrea.calogero@unimib.it)} and
R. Pini\thanks{Dipartimento di Metodi Quantitativi per le Scienze
Economiche ed Aziendali, Universit\`a degli Studi di Milano
Bicocca, Via Bicocca degli Arcimboldi 8, I-20126 Milano
\tt(rita.pini@unimib.it)}}

\title{$c$ horizontal convexity on Carnot groups}

\maketitle
\begin{abstract}
\noindent Given a real--valued function $c$ defined on the
cartesian product of a generic Carnot group $\G$ and the first
layer $V_1$ of its Lie algebra, we introduce a notion of $c$
horizontal convex ($c$ H--convex) function on $\G$ as the supremum
of a suitable family of affine functions; this family is defined
pointwisely, and depends strictly on the horizontal structure of
the group. This abstract approach provides $c$ H--convex functions
that, under appropriate assumptions on $c,$ are characterized by
the nonemptiness of the $c$ H--subdifferential and, above all, are
locally H--semiconvex, thereby admitting horizontal derivatives
almost everywhere. It is noteworthy that such functions can be
recovered via a Rockafellar technique, starting from a suitable
notion of $c$ H--cyclic monotonicity for maps. In the particular
case where $c(g,v)=\langle \xi_1(g),v \rangle,$ we obtain the
well--known weakly H--convex functions introduced by Danielli,
Garofalo and Nhieu. Finally, we suggest a possible application to
optimal mass transportation.
\end{abstract}


\section{Introduction}

In $\R^n$ and, more generally, in a Banach space $\mathbf{X}$, the
notion of convexity of a function $f:\mathbf{X}\to
(-\infty,+\infty]$ can be given in terms of the supremum of the
affine functions $x\mapsto \langle x,y\rangle +\alpha$ lying below
the function itself. Among the nice properties enjoyed by proper,
lower semicontinuous and convex functions, we recall that they can
be characterized by the nonemptiness of the subdifferential
$\partial f$ at every point of their domain and, by a well--known
result due to Rockafellar, they can be completely recovered by
their subgradients. In addition, the multivalued map $x\mapsto
\partial f(x)$ benefits from an interesting condition, since
it can be characterized as a maximal monotone map \cite{Ro1970}.

This abstract formulation of convexity is fit for an extension by
substituting the affine function $\langle \cdot, y\rangle +\alpha$
with the more general function $c(\cdot,y)+\alpha,$ where
$c:\mathbf{X}\times \mathbf{X}^*\to\R.$ This generalization leads
to the definition of $c$ convex function, and the associated $c$
subdifferential multivalued map $x\mapsto
\partial_cf(x)$ arises in a natural way, together with the notion of $c$
cyclic monotonicity. These concepts date back to a first paper by
E.J. Balder \cite{Ba1977}, and were introduced in order to extend
the duality theory to nonconvex optimization problems. In this
framework, H. Dietrich \cite{Di1988} investigated several properties
of $c$ subdifferentiability and local $c$ subdifferentiability of
$c$ convex functions. Subsequently, L. R\"uschendorf \cite{Ru1996},
in connection with the coupling problem, gave a characterization of
optimal solutions via generalized subgradients of $c$ convex
functions.

As a matter of fact, as can be found in the fundamental paper by
W. Gangbo and R.J. McCann \cite{GaMc1996}, a context where $c$
concavity plays a key role is in finding the solution of an
optimal mass transportation problem, where $c(x,y)$ denotes the
cost per unit mass displaced from $x$ to $y$ ($x,y\in\R^d$).
Indeed, the support of the optimal measure on $\R^d\times\R^d$ is
contained in the graph of $x\mapsto
\partial^c \psi(x),$ where $\psi$ is a $c$ concave function called potential; if $c(x,y)=h(x-y),$ then,
under suitable regularity assumptions on $h,$ like strict
convexity and superlinearity, a deep result says that this
multivalued map $\partial^c \psi$ is essentially single--valued. A
remarkable result concerns the finiteness and the regularity of a
$c$ concave function; it is noteworthy that a $c$ concave function
inherits structure and smoothness from the function $c,$ like
locally Lipschitz, local semiconcavity, and local boundedness.
This implies that there exists an optimal transport map $s$
defined on $\mathrm{dom}(\nabla\psi)\subset\R^d$ by the formula
$s(x)=x-(\nabla h)^{-1}(\nabla \psi(x)).$

In the quite recent literature, in the Heisenberg group $\H$ and,
more generally, in Carnot groups, several concepts of convexity
have been introduced (see, for instance, \cite{DaGaNh2003}, or
\cite{CaDaPaTy2007}). Among them, the most suitable to many
purposes is the so--called weakly H--convexity (H--convexity, in
the sequel). An H--convex function $u$ is, essentially, a function
that is convex along any horizontal line, a particular horizontal
curve. Balogh and Rickly proved that these functions are regular
enough, since they are locally Lipschitz continuous with respect
to any homogeneous distance (see \cite{BaRi2003}, \cite{Ri2006}).
In \cite{CaPi2008}, we show that, for real--valued functions on
$\H$, H--convexity is equivalent to H--subdifferentiability, i.e.
the horizontal subdifferential is nonempty at every point of the
domain; in this paper this result will be extended to a generic
Carnot group.

Unexpectedly, it turns out (see \cite{CaPi2010}) that there is an
abstract definition in $\H$ of convexity, given in terms of
H--affine functions, that is equivalent to the H--convexity. A
real--valued function $u$ is abstract H--convex if
$$
u(g)=\sup_{(v,\alpha)\in \Pca_g} (\langle \xi_1(g),v\rangle+\alpha)
$$
where $\xi_1$ is defined in Section \ref{basic},
 and $\Pca_g$ is
the set of pairs $(v,\alpha)\in V_1\times\R$ such that the
H--affine function $g'\mapsto \langle \xi_1(g'),v\rangle+\alpha$
supports $u$ on every horizontal line through $g$ (see
\cite{DaGaNh2003}, p. 320). Let us stress the peculiarity of the
set of parameters $(v,\alpha)\in V_1\times\R,$ that depends on the
point $g.$ Moreover, let us notice that $v$ belongs to the first
layer of the Lie algebra of $\G$ and plays the role of a \lq\lq
subgradient".

This point of view can be extended, by taking a general function
$c(g,v)$ instead of $\langle \xi_1(g),v\rangle;$ this paper is
devoted to the study of the main features of these functions, that
will be called $c$ horizontal convex ($c$ H--convex, briefly), as
well as to the interrelationships with their $c$
H--subdifferentials. We cannot leave unmentioned the papers by L.
Ambrosio and S. Rigot \cite{AmRi2004}, \cite{Rig2005}, where the
optimal transport mass is investigated in the framework of
particular Carnot groups; in this context they introduce a
different notion of $c$ concavity and $c$ superdifferential that
does not take into account the \lq\lq horizontal" structure of a
Carnot group.

Our investigation follows some classical steps in convex analysis.
The main results mirror similar ones in the context of classical
convexity. First of all, in Section \ref{section H-convexity}, we
provide the mentioned link between H--convexity of a real--valued
function and H--subdifferentiability in a generic Carnot group
$\G$ (see Theorem \ref{converse}). In Section \ref{section about c
H-convexity}, we introduce the notions of $c$ H--convexity and $c$
H--subdifferential $\partial_H^c u$ for a proper function $u:\G\to
(-\infty, +\infty]$. Theorem \ref{converse c} characterizes a
real--valued $c$ H--convex function via its $c$
H--subdifferentiability; this result is not an extension of
Theorem \ref{converse}, since it deals with abstract H--convexity.

Section \ref{section regularity} is devoted to the problem of the
regularity of a proper $c$ H--convex function, and we try to
establish the almost sure single--valuedness of its $c$
H--subdifferential. In Euclidean spaces, semiconvexity turns out
to be a fundamental tool for the study of $c$ convex functions
(see, for instance, \cite{Vil2009}, Chapter 10). Semiconvexity can
be extended in a natural way in a Carnot group starting from
H--convexity, and it gives rise to the notion H--semiconvexity
(see Definition \ref{locally H semiconvex}). Despite their
abstract and entangled definition, $c$ H--convex functions prove
to be well--behaved whenever $c$ is. One of the most interesting
result of the paper is Theorem \ref{regolarità semiconvex}: we
show that, like in the classical case, our functions are locally
H--semiconvex, and therefore they share the regularity of the
H--convex functions. This entails that, in the real--valued case,
and under measurability assumptions if the step is greater than 2,
a $c$ H--convex function $u$ is differentiable almost everywhere
along the horizontal directions; furthermore, we get that
$\partial_H^c u(g)$ is a singleton for almost every $g.$

Another relevant issue that shares its aim with classical
convexity concerns the connection between $c$ H--convexity of a
function on $\G,$ and the $c$ H--cyclic monotonicity of a subset
of $\G\times V_1$ (see Definition \ref{definizione c Hcyclically
monotone}). The main and more delicate outcome of Section
\ref{section about c H-convexity} shows that, from a $c$
H--cyclically monotone set and via Rockafellar techniques, it is
possible to construct, at least locally, a $c$ H--convex function
$u$ such that the graph of $g\mapsto\partial_H^c u(g)$ contains
the starting set. In this setting, the analysis of the finiteness
of the function plays a critical role.

Finally, inspired by the precious paper by Gangbo and McCann
\cite{GaMc1996}, we present a possible application of all these
arguments and tools to the optimal mass transportation problem in
the Heisenberg group. Despite this application arises in a very
particular situation, where the optimal map moves the points only
along horizontal segments, we think that our approach could be
potentially interesting.

\section{Basic notions on Carnot groups}\label{basic}

 A Carnot group $\mathbf{G}$ of step $r$ is a connected, simply
 connected, nilpotent
Lie group whose Lie algebra $\frak{g}$ of left--invariant vector
fields admits a stratification, i.e. there exist non zero subspaces
$\{V_j\}_1^r$ such that
\begin{itemize}
\item[] $\frak{g}=V_1\oplus V_2\oplus\ldots\oplus V_r,$

\item[] $ [V_1,V_j]=V_{j+1}\qquad\qquad j=1,\ldots r-1,$

\item[]
$[V_1,V_r]=0.$
\end{itemize}

We assume that a scalar product $\langle\cdot,\cdot\rangle$ is
given on $\frak{g}$ for which the levels $V_j$ are mutually
orthogonal. The first layer $V_1$ of the Lie algebra plays a key
role: we call \emph{horizontal vector fields} its elements, and
denote by $m$ its dimension.

\noindent We fix an orthonormal basis $X=\{X_1,X_2,\ldots,X_m\}$
of $V_1,$ and we continue to denote by $X$ the corresponding
system of left--invariant vector fields on $\G$ defined by
$X_i(g)=(L_g)_* (X_i),\ i=1,\ldots,m,$ where $(L_g)_*$ is the
differential of the left translation on $\G$ defined by
$L_g(g')=gg'$. The system $X$ defines a basis for the horizontal
sub--bundle $\Hca\G$ of the tangent bundle $\Tca\G$ (i.e. $\Hca_g=
{\rm span}\{X_1(g),\ldots,X_m(g)\}$ for every $g\in\G$).

\noindent The action of $X_i$ on a function $u:\mathbf{G}\to
\mathbb{R}$ is given by
$$
X_i u(g)=\lim_{\alpha\to 0}\frac{u(g\, \texttt{\rm exp}(\alpha
X_i))-u(g)}{\alpha}.
$$
Clearly,  $\exp: \frak{g}\to \mathbf{G}$ is the exponential map,
 a global diffeomorphism; we denote by
$\xi=(\xi_1,\xi_2,\ldots,\xi_r)$ the inverse of $\exp,$ where
$\xi_j:\G\to V_j.$

A natural family of non--isotropic dilations on $\frak{g}$
associated with its grading is given by
$\Delta_\lambda(v_1+v_2+\ldots+v_r)=\l v_1+\l^2 v_2+\ldots+\l^r
v_r,$ if $v_i\in V_i$ for every $1\le i\le r.$ By means of the
exponential map, one lifts these dilations to the family of the
automorphisms $\delta_\l(g)=\exp( \Delta_\lambda( \xi(g))).$
 The homogeneous dimension associated with the dilations
$\{\delta_\l \}_{\l>0}$ is given by $Q=\sum_{i=1}^r i \,\texttt{\rm
dim} V_i$ that often replaces the topological dimension
$N=\sum_{i=1}^r \,\texttt{\rm dim} V_i$ in the study of Carnot
groups.

The Euclidean distance to the origin $|\cdot|_\frak{g}$ on
$\frak{g}$ induces a homogeneous pseudo--norm $\|\cdot\|_\frak{g}$
on $\frak{g}$ defined by $\|v_1+v_2+\ldots +
v_r\|_\frak{g}=\left(\sum_{i=1}^r
|v_i|_\frak{g}^{2r!/i}\right)^{2r!}.$ Again, via the exponential
map, we lift $\|\cdot\|_\frak{g}$ to a pseudo--norm
$\|\cdot\|_\G,$ and hence to a pseudo--distance $d$ on $\G$
defining $\|g\|_\G=\|\xi(g)\|_\frak{g}$ and $d(g,g')=\|g^{-1}
g'\|_\G.$

Let $\Omega\subset\G$ be an open set, $k$ be a non negative
integer, and $0<\alpha\le 1.$ The class $\Gamma^k(\Omega)$
represents the Folland--Stein space of functions having continuous
derivatives up to the order $k$ with respect to the horizontal
vector fields $X_1,\ldots, X_m.$ A function $u:\Omega\to\R$ is
said to belong to the class $\Gamma^{0,\alpha}(\Omega)$ if there
exists a positive constant $C_\alpha$ such that $$|u(g)-u(g')|\le
C_\alpha d(g,g')^\alpha,$$ for every $g$ and $g'$ in $\Omega.$ A
function $f\in\Gamma^{1}(\Omega)$ belongs to the class
$\Gamma^{1,\alpha}(\Omega),$ if for every $i=1,\ldots,m,$ the
horizontal derivative $X_i f$ exists in $\Omega$ and
 $X_i f\in\Gamma^{0,\alpha}(\Omega).$
As usual, we say that $u$ is Lipschitz continuous if $u\in
\Gamma^{0,1};$ the symbol $\Gamma^{0,1}_{\rm loc}(\Omega)$ denotes
the class of locally Lipschitz continuous functions on $\Omega.$

Let us recall that the horizontal gradient of a function
$u\in\Gamma^1(\Omega)$ at $g\in\Omega$ is the element of $V_1$
$$\X u(g)=\sum_{i=1}^m (X_iu(g))X_i.$$

The \emph{horizontal plane} $H_{g}$ associated to $g\in \G$ is
given by
\begin{equation}\label{piano orizzontale}
H_{g}=L_{g}\left(\exp(V_1)\right)=\{ g'\in\G: g'=g h,\ \texttt{\rm
with}\ h\in\exp(V_1)\}.
\end{equation}
Note that $g'\in H_{g}$ implies that $g\in H_{g'}$ and
$g^{-1}g'\in H_{e},$ where $e$ is the unit element of the group
$\G. $ If we consider the set $H_e,$ and identify $\G$ with $\R^N$
(remember that $N$ is the topological dimension), it turns out
that the set $H_e$ is an iperplane in $\R^N.$ Differently, if
$g\neq e,$ one can show that the horizontal plane $H_{g}$ is an
iperplane in the classical sense (in particular an $\R^N$--convex
set, using the subsequent notation) if and only if $\G$ has step 2
(see Example \ref{engel}).

As a matter of fact, the elements of the first layer $V_1$ of the
Lie algebra $\frak{g}$ generate all the vector fields of
$\frak{g}$ and consequently, via the exponential map, the points
of the horizontal plane $H_e$ play a similar role in $\G.$ More
precisely, the following structure result holds:

\begin{proposition}\label{FoSt}(see \cite{FoSt1982}, Lemma 1.40).
Let $\G$ be a stratified group. Then, there exist $C>0$ and
$R\in \mathbb{N}$ such that any $g\in \G$ can be expressed as
$g=h_1h_2\cdots h_R,$ with suitable $h_i\in H_e$ and
$\|h_i\|_{\G}\le C\|g\|_\G,$ for every $i=1,2,\dots, R.$
\end{proposition}

We recall that a Lipschitz curve $\gamma:[0,T]\to \G$ is said to be
 horizontal if $\gamma'(\l)\in\Hca_{\gamma(\l)},$ i.e. $
 \gamma'(\l)=\sum_{i=1}^m a_i(\l)X_{\gamma(\l)},$
  for almost every
$\l\in[0,T].$ The sub--Riemannian length of a horizontal curve
$\gamma$ is $$L(\gamma)=\int_0^T \left(\sum_{i=1}^m
a_i^2(\l)\right)^{1/2} d\l;$$ the Carnot--Caratheodory distance
$d_{CC}$ from $g$ to $g'$ is
$$
d_{CC}(g,g')=\inf\{L(\g): \g\ \texttt{\rm is a horizontal curve
connecting}\ g \  \texttt{\rm and}\ g'\}.
$$
A curve $\g$ joining $g$ and  $g'$ is a geodesic if it is a length
minimizing horizontal curve, i.e. $L(\g)=d_{CC}(g,g').$ Another
kind of curve connecting two points $g$ and $g'$ arises as their
twisted convex combination $\sigma_{g,g'}$ defined by
\begin{equation}\label{twiset convex combination}
\sigma_{g,g'}(\l)=g\delta_{\lambda}(g^{-1}g'),\quad \lambda\in
[0,1].
\end{equation}
If $g'\in H_g,$ we say that $\sigma_{g,g'}$ is a \emph{horizontal
segment}; it is a horizontal curve and, in particular, a geodesic.

We say that  $A\subset \R^n$ is $\R^n$--convex if $(1-\l) x+\l
y\in \Omega,$ for every $x,y$ in $\Omega,$ and $\l\in[0,1].$
Consequently, a function is $\R^n$--convex if $u((1-\l) x+\l y)\le
(1-\l)u(x)+\l u(y),$ with $x,y,\l$ as before. An $\R^n$--segment
is the $\R^n$--convex hull of two points, and an $\R^n$--plane is
the set $\{x\in\R^n:\ \langle x,a\rangle={b}\},$ for some fixed
${a}\in\R^n$ and ${b}\in\R.$ These notations should be pedant, but
it is important in this paper to distinguish the different notions
of convexity, plane, segment that we introduce.

Let us explain these arguments with two basic models.

\footnotesize
\begin{example}\label{heisenberg}\rm The Heisenberg group $\H.$

\noindent The Heisenberg group $\H$ is the Lie group whose
 Lie algebra $\frak{h}$ admits a stratification of step 2; in particular
 $\frak{h}=\mathbb{R}^3=V_1\oplus V_2,$
with
\begin{equation}\label{algebra Heisenberg}
\begin{array}{ll}
V_1={\rm span}\left\{X_1,X_2\right\}\quad&\texttt{\rm with}\ \
X_1=\partial_{x}-\frac{y}{2}\partial_{t}\quad \texttt{\rm and}\ \
X_2=\partial_{y}+\frac{x}{2}\partial_{t},\\
  V_2={\rm
span}\left\{T\right\}&\texttt{\rm with}\ \ T=\partial_{t}.
\end{array}
\end{equation}
The bracket $[\cdot ,\cdot ]:\frak{h} \times\frak{h} \to \frak{h}$
is defined as $[X_1,X_2]=T,$ and it vanishes in the other cases;
taking into account the action of the bracket, $X*Y$ is defined by
the Baker--Campbell--Dynkin--Hausdorff formula
\begin{equation}\label{baker campbell dynkin hausdorff}
X*Y=X+Y+[X,Y]/2.
\end{equation}
 The exponential map $\exp(\alpha X+\beta Y+\gamma T)=
(\alpha,\beta,\gamma)$ enjoys the property
$\exp(X)\exp(Y)=\exp(X*Y),$ for every $X$ and $Y$ in $\frak{g};$
consequently, the law group on $\H$ is
$$
gg'=(x,y,t)(x',y',t')= (x+x',y+y',t+t'+(xy'-x'y)/2).$$ The dilation
is a family of automorphisms given by
$\delta_\lambda(x,y,t)=(\lambda x,\lambda y,\lambda^2 t),$ and hence
the homogeneous dimension is 4. Given two points $g=(x,y,t)$ and
$g'=(x',y',t'),$ the non commutative twisted convex combination in
(\ref{twiset convex combination}) is
$$
\sigma_{g,g'}(\l)= \left( (1-\lambda)x+\lambda x',
(1-\lambda)y+\lambda y',
t+\lambda(xy'-x'y)/2+\lambda^2(t'-t+(x'y-xy')/2) \right).
$$
The horizontal plane $H_g$ is, by (\ref{piano orizzontale}),
$$
H_{g}=\left\{(x',y',t')\in \H:\  t'=t+(xy'-x'y)/2,\ \ x',y'\in
\mathbb{R}\right\};
$$
it is a \lq\lq real\rq\rq\ plane, i.e. an $\R^3$--plane. If we
choose $g'$ on the horizontal plane $H_g,$ the curve
$\sigma_{g,g'}$ is a horizontal curve and a geodesic that we call,
by definition, horizontal segment $\gamma$ from $g$ to $g':$ more
precisely,
$$
\gamma(\lambda) =( (1-\lambda)x+\lambda x', (1-\lambda)y+\lambda y',
t+\lambda(xy'-x'y)/2) .
$$
Note that $\gamma$ is an $\mathbb{R}^3$--segment lying in $
H_{g}\cap H_{g'}.$
\end{example}

\begin{example}\label{engel}\rm The Engel group $\E.$

 \noindent The Engel group is a Carnot group of step 3 and, in
some sense, is an extension of $\H:$ indeed if we consider the Lie
algebra $\frak{e}=\mathbb{R}^4=V_1\oplus V_2\oplus V_3$ defined
by, using (\ref{algebra Heisenberg}),
$$
\begin{array}{ll}
V_1={\rm span}\{\tilde{X}_1,\tilde{X}_2\}\quad&\texttt{\rm with}\ \
\tilde{X}_1=X_1 -(\frac{t}{2}+\frac{xy}{12})\partial_s\quad
\texttt{\rm and}\ \
\tilde{X}_2=X_2 +\frac{x^2}{12}\partial_s,\\
  V_2={\rm
span}\{\tilde{T}\}&\texttt{\rm with}\ \
\tilde{T}=T+\frac{x}{2}\partial_s,\\
 V_3={\rm
span}\{\tilde{S}\}&\texttt{\rm with}\ \ \tilde{S}=\partial_s.
\end{array}
$$
The bracket acts as $[\tilde{X}_1,\tilde{X}_2]=\tilde{T},\ \
[\tilde{X}_1,\tilde{T}]=\tilde{S},$ and it vanishes in the other
cases. Since, in $\frak{e},$  in the
Baker--Campbell--Dynkin--Hausdorff formula (\ref{baker campbell
dynkin hausdorff}) there is one more term (precisely
$([X,[X,Y]]+[Y,[Y,X]])/12$) and $\texttt{\rm exp}(\alpha
\tilde{X}+\beta \tilde{Y}+\gamma \tilde{T}+\eta
\tilde{S})=(\alpha,\beta,\gamma,\eta),$ the  group law in $\E$
becomes
 $$gg'=
\Bigl(x+x',y+y',t+t'+(xy'-x'y)/2,
s+s'+(xt'-x't)/2+(x-x')(xy'-yx')/12 \Bigr),$$ where $g=(x,y,t,s)$
and $g'=(x',y',t',s').$ The horizontal plane $H_{g}$ is
$$
\begin{array}{ll}
\quad H_{g}=\Bigl\{(x',y',t',s')\in \E:\quad   t'=t+(xy'-x'y)/2,\\
\qquad\qquad\qquad s'=s+(-6t(x'-x)+2x^2y'-2x'xy+yx'^2-xx'y')/12, \
\texttt{\rm with}\ x,y\in \mathbb{R} \Bigr\};
\end{array}
$$
note that $H_{g}$ is not an $\mathbb{R}^4$--plane. Clearly, the
dilation is given by
 $\delta_\lambda(x,y,t,s)=(\lambda x,\lambda y,\lambda^2
 t,\lambda^3s).$
If we consider $g'\in H_{g},$ the horizontal segment $\gamma$ with
endpoints $g$ and $g'$ is defined via (\ref{twiset convex
combination}): $\gamma$ is a geodesic, lies in $H_g\cap H_{g'}$
and, in general, is not an $\mathbb{R}^4$--segment.
\end{example}
\normalsize

We have seen that, unlike in the Euclidean spaces, where the
Euclidean distance is the most natural choice, in a Carnot group
several distances were introduced for different purposes. However,
all of these distances $\rho$ are homogeneous, namely, they are
left invariant and satisfy the relation $\rho(\delta_r g',\delta_r
g)=r\rho(g',g)$ for every $g',g\in \G,$ and $r>0.$ The distance
functions $d$ and $d_{CC}$ are homogeneous, equivalent, and have
the same value at the endpoints of a horizontal segment.

Let $\rho$ be any homogeneous distance on $\G,$ and let $u:\G\to
\R.$ We say that $u$ is Pansu differentiable at $g\in\G$ if there
exists a $\G$--linear map $Du(g):\G\to\R,$ i.e., a group
homomorphism that satisfies the relation $Du(g)(\delta_r h)=r
Du(g)(h)$ for every $h\in\G$ and $r>0,$ and
$$ \lim_{\rho(h,e)\to
0}\frac{|u(gh)-u(g)-Du(g)(h)|}{\rho(h,e)}=0.
$$
We call the map $Du(g)$ the Pansu differential of $u$ at $g.$
 An easy computation gives us that if $u$ is Pansu differentiable at
 $g$, then
 $$
 Du(g)(h)=\lim_{\lambda\to
0^+}\frac{u(g\delta_\lambda(h))-u(g)}{\lambda} $$
exists for every
$h\in\G.$ If $u\in \Gamma^1(\G),$ then the Pansu differential
$Du(g)$ is given by the formula
$$
Du(g)(h)=\langle\X u(g),\xi_1(h)  \rangle,
$$
for every $g$ and $h$ in $\G$ (see \cite{DaGaNh2003}).

It is known that a Rademacher--Stefanov type result holds in the
Carnot group setting; therefore, a Lipschitz continuous function
is differentiable almost everywhere in the horizontal directions.
A further result, due to Danielli, Garofalo and Salsa, will play a
crucial role in the sequel:
\begin{theorem}\label{DaGaSa}(see \cite{DaGaSa2003}, Theorem 2.7).
Let $\Omega$ be an open subset of $\G,$ and $u:\Omega\to \R,$ with
$u\in \Gamma^{0,1}(\Omega).$ Then there exists a set $E\subset
\Omega$ with Haar measure zero such that the Pansu differential
$Du(g)$ and the horizontal gradient $\X u(g)$ exist for every $g\in
\Omega\setminus E,$ and
$$
Du(g)(h)=\langle \X u(g),\xi_1(h)\rangle,\qquad \texttt{\it for
every}\ h\in \G. $$
Furthermore, $\X u\in L^{\infty}(\Omega).$
\end{theorem}

Finally, for what concerns classical convex analysis, we will
refer to \cite{Ro1969}. In particular, we say that a function $u$
defined on a subset $\Omega$ of $\G$ is \emph{proper} if $u(g)\neq
-\infty$ for every $g\in \Omega,$ and $u\not\equiv +\infty;$
moreover, if $u(g)\neq +\infty$ for every $g\in\Omega,$ then we
say that $u$ is \emph{real--valued}. The domain of $u,$
$\texttt{\rm dom}(u),$  is the subset of $\Omega$ where $u$ is
finite.

\section{H--convexity and H--subdifferentiability}\label{section
H-convexity}

In the last few years, several notions of convexity have been
introduced in the framework of Carnot groups, but the most
suitable one showed to be the notion of H--convexity. This notion
is due to Caffarelli in unpublished works from 1996, and it
appeared in the paper \cite{DaGaNh2003}; afterwards, several
papers have been devoted to the investigations of H--convexity.
Among other things, the horizontal Monge--Amp\`ere equation in
$\G$ defined by $\det [\X^2 u]^*(g)=f(g,u,\X u)$ is (degenerate)
elliptic precisely on the class of $u\in \Gamma^2(\G)$ which are
H--convex. This section will concern results about H--convex
functions in a Carnot group $\G.$

A subset $\Omega$ of $\mathbf{G}$ is {\it H--convex} if it
contains every horizontal segment with endpoints in $\Omega,$ i.e.
$ g\delta_{\l}(g^{-1}g')\in \Omega,$ for every $g\in \Omega,\
g'\in H_g\cap \Omega$ and $\l\in [0,1].$
\begin{definition} {\rm Let $\Omega\subset \mathbf{G}$ be H--convex.
A function $u:\Omega\to (-\infty,+\infty]$ is {\it H--convex} if
it is $\R$--convex on every horizontal segment, i.e.
\begin{equation}\label{def Hconvex}
u(g\delta_{\l}(g^{-1}g'))\le (1-\lambda)u(g)+\lambda u(g')
\end{equation}
 for all $g\in\Omega$, $g'\in H_g\cap\Omega,$ and
$\lambda\in[0,1].$}
\end{definition}
A function $u:\Omega\subset \G\to[-\infty,+\infty)$ is said to be
H--concave if $-u$ is H--convex.

It is clear, by the definition, that an $\R^3$--convex function
$u$ is H--convex in the Heisenberg group, since every horizontal
segment is a particular $\R^3$--segment. This argument can be
extended to any Carnot group of step two. On the contrary, if one
consider a group $\G$ of step greater than 2 this is no longer
true. An enlightening example can be found in \cite{Ri2006}: the
function $u:\E\to \R,$ $u(x, y,t,s)=s$ is not H--convex in the
Engel group $\E$, despite it is $\R^4$--convex.

In spite of the notion of H--convexity, that requires a suitable
behaviour on the horizontal planes only, H--convex functions enjoy
some nice regularity properties. Balogh and Rickly (see
\cite{BaRi2003} if $\G=\H,$ and \cite{Ri2006}) proved the
following result:

\begin{theorem}\label{Ri} (see \cite{Ri2006}, Theorem 1.4).
Let $\Omega\subset\mathbf{G}$ be an H--convex, open subset. Then
every H--convex function $u:\Omega\to \R,$ measurable if the step
of $\mathbf{G}$ is greater than 2, belongs to
$\Gamma^{0,1}_{\mathrm{loc}}(\Omega).$
\end{theorem}
The possibility to remove the measurability assumption in the
previous result, is an interesting and open question.

In \cite{DaGaNh2003}, a regular function $u:\Omega\to\R,$ where
$\Omega$ is an open and H--convex subset of $\G,$ is characterized
in terms of its horizontal gradient $\X u$, and its symmetrized
horizontal Hessian $[\X^2 u]^*.$ Indeed, if $u\in
\Gamma^1(\Omega),$ then $u$ is an H--convex function if and only
if
\begin{equation}\label{condizione primo ordine}
u(g')\ge u(g)+\langle \X u(g),
\xi_1(g')-\xi_1(g)\rangle,\qquad\forall g\in\Omega,\quad \forall
g'\in H_g\cap\Omega;
\end{equation}
if $u\in \Gamma^2(\Omega),$ then $u$ is H--convex if and only if
$[\X^2 u]^*(g)$ is positive semidefinite for every $g\in \Omega,$
where
$$
[\X^2 u]^*(g)=\frac{1}{2}\{\X^2u(g)+\X^2u(g)^T\}
$$
and $\X^2u(g)=(X_iX_ju(g))_{i,j=1,\ldots, m}$ is an $m\times m$
matrix.

It is well known that if $f:\R^n\to\R$ is a differentiable
function, then the $\R^n$--convexity of $f$ can be characterized
by the monotonicity of the gradient, i.e., $(\nabla f(x)-\nabla
f(y)) \cdot (x-y)\ge 0,$ for every $x,y$ in the domain (see, for
instance, \cite{AvrDieSchZa1988}, Theorem 2.13). This result can
be adapted to the sub--Riemannian setting; indeed, if
$u\in\Gamma^1(\Omega),$ then one can easily show from
(\ref{condizione primo ordine}) that $u$ is an H--convex function
if and only if
\begin{equation}\label{caratterizzazione con gradiente H}
\langle \X u(g)-\X u(g'), \xi_1(g)-\xi_1(g')\rangle \ge
0,\qquad\forall g\in \Omega,\quad \forall g'\in H_g\cap\Omega.
\end{equation}
We say that the set $\{g_i\}_{i=0}^n\subset \Omega$ is an
\emph{H--sequence} if, for some $n>0$ and for every $i=0,\dots,
n-1,$ $g_{i+1}\in H_{g_i}.$  An H--sequence is closed if $g_n\in
H_{g_0};$ in this case, we usually set $g_{n+1}=g_0.$
 This
notion, that will be fundamental in the next sections, allow us to
extend the characterization in (\ref{caratterizzazione con
gradiente H}). Indeed, an easy calculation shows that, when
$u\in\Gamma^1(\Omega),$ then $u$ is an H--convex function if and
only if
\begin{equation}\label{caratterizzazione con monotomia gradiente H}
\sum_{i=0}^{n} \langle \X u(g_i), \xi_1(g_{i+1})\rangle \le
\sum_{i=0}^{n} \langle \X u(g_i), \xi_1(g_{i})\rangle,
\end{equation}
for every closed H--sequence $\{g_i\}_{i=0}^n\subset\Omega.$ This
last property will lead to consider the more general notion of $c$
H--cyclic monotonicity in Section \ref{$c$ H--cyclic
monotonicity}.

In \cite{DaGaNh2003} the authors relate the property of H--convexity
of a real--valued function to the nonemptyness of its
H--subdifferential. Let us recall that the H--subdifferential of a
function $u:\Omega\subset\mathbf{G}\to(-\infty,+\infty]$ at $g_0\in
\Omega$ is defined as
$$
\partial_Hu(g_0)=\{p\in V_1:\,u(g)\ge u(g_0)+\langle
p,\xi_1(g)-\xi_1(g_0)\rangle,\,\forall g\in H_{g_0}\cap\Omega\}.
$$
Moreover, we say that $\partial^Hu(g_0)$ is the
H--superdifferential of $u$ at $g_0$ if $\partial^Hu(g_0)=
-\partial_H(-u)(g_0).$

A first link between H--subdifferentiability of a function and
H--convexity is provided by the following:
\begin{proposition}
\label{prop0} (see \cite{DaGaNh2003},  Proposition 10.5). Let
$u:\Omega\to\R,$ where $\Omega$ is an open and H--convex subset of
$\mathbf{G}.$ If $\partial_Hu(g)\neq \emptyset$ for every $g\in
\Omega,$ then $u$ is H--convex.
\end{proposition}
The converse of this result, as in the classical case, is more
difficult.
 In \cite{CaPi2008} we prove that this
holds when $\mathbf{G}=\H$. As a matter of fact, next theorem shows
that the result can be improved.

\begin{theorem}\label{converse}
Let $u:\Omega\subset \mathbf{G}\to\R,$ where $\Omega$ is open and
H--convex. Let $u$ be H--convex, and measurable if $r>2.$ Then
$\partial_Hu(g)\neq \emptyset$ for every $g\in \Omega.$
\end{theorem}
In order to prove  Theorem \ref{converse}, one can extend in a
natural way the proof in \cite{CaPi2008} from the Heisenberg group
to a generic Carnot group, with the additional assumption that $u$
is measurable if $r>2.$ We recall the main tools of such proof and
we leave its details to the reader. In this setting, a fundamental
role is played by the regularity results for H--convex functions
due to Balogh and Rickly (see Theorem \ref{Ri}), and the
differentiability almost everywhere in the horizontal directions
for Lipschitz continuous functions due to Danielli, Garofalo and
Salsa (see Theorem \ref{DaGaSa}). The assumption of H--convexity
and the two results above, together, lead to the inclusion $\X
u(g)\in
\partial_Hu(g)$ a.e. in $\Omega.$
A crucial point lies in proving that the graph of the multivalued
map $g\mapsto
\partial_Hu(g)$ is closed, i.e.
for every sequence $\{(g_n,p_n)\}_n\subset \Omega\times V_1$ with
$p_n\in
\partial_Hu(g_n),$ such that
$ g_n\to g_0\in \Omega$ and $p_n\to p,$ then $p\in
\partial_Hu(g_0).$
In order to do this, we exploit the continuity of the function
$u$, together with the \lq\lq continuity" of the left translation
on the group, that is involved in the definition of the horizontal
planes (\ref{piano orizzontale}); more precisely, given $g_0\in
\Omega,$ and $g'\in H_{g_0}\cap \Omega,$ for every $g_n\to g_0$
there exists $\{g'_n\}_n$ such that $g'_n\in H_{g_n}\cap\Omega$
and $g'_n\to g'.$ The reader can give a look at Lemma 4.1 in
\cite{CaPi2008} to find more details in the case $\G=\H.$

\section{ \emph{c} H--convexity and \emph{c}
H--subdifferential}\label{section about c H-convexity}

The class of $c$ convex functions was introduced, to our
knowledge, by Dietrich \cite{Di1988}, and subsequently exploited
by several authors in connection with optimal couplings and
optimal mass transportation problems; to get an idea about it, one
can read the paper by R\"uschendorf \cite{Ru1996}, or give a look
at the book by C. Villani \cite{Vil2009}. Briefly, if
$\Omega_1,\Omega_2$ are two sets, and $c:\Omega_1\times
\Omega_2\to \R,$ then a proper function $f:\Omega_1\to
(-\infty,+\infty]$ is said to be $c$ convex if there exists a set
$\Pca\subset \Omega_2\times \R$ such that
\begin{equation}\label{c convexity generale}
f(x)=\sup_{(y,\alpha)\in\Pca}(c(x,y)+\alpha),\qquad\forall x\in
\Omega_1.
\end{equation}
In the investigation about the properties of $c$ convex functions,
a fundamental role is played by the notion of $c$ subdifferential
$\partial_c f$ defined as
\begin{equation}\label{c subdifferential generale}
\partial_c f(x)=\{y\in \Omega_2:\ f(x')\ge f(x)+c(x',y)-c(x,y),\ \forall x'\in
\Omega_1\}. \end{equation} In particular cases, for instance if
$\Omega_1=\Omega_2=\R^n$ and $c(x,y)=\langle x,y\rangle,$ one can
easily recover some classical notions: in (\ref{c convexity
generale}) we obtain the abstract notion of convexity, where a
convex function is defined as the pointwise supremum of a family
of affine functions; in (\ref{c subdifferential generale}) we
obtain the notion of subgradient, i.e. the set of coefficients $y$
such that the affine function $x'\mapsto \langle
x'-x,y\rangle+f(x)$ supports the function $f$ at the point $x.$

In \cite{AmRi2004}, the authors deal with an optimal mass
transportation problem in the Heisenberg group, and they are lead
to consider the class of $c$ convex functions on $\H.$  In
particular they prove the existence and the uniqueness of an
optimal transport map assuming that the cost function
$c:\H\times\H\to\R$ is either the function $d^2,$ or the function
$d^2_{CC}$ (see \cite{Rig2005} for the more general case of groups
of type $H$). As a matter of fact, the notion of $c$ convexity
they work with does not take into account the horizontal
structure; more precisely, they say that $f:\H\times\H\to\R$ is
$c$ convex if (\ref{c convexity generale}) holds, at every
$x\in\H,$ for a suitable nonempty set $\Pca\subset\H\times\R.$
Consequently, their definition of $c$ subdifferential is exactly
as in (\ref{c subdifferential generale}), with
$\Omega_1=\Omega_2=\H.$

The aim of this paper is the investigation of $c$ convexity from
another viewpoint: in Sections \ref{section about c
H-convexity}--\ref{$c$ H--cyclic monotonicity} we provide a
different notion of $c$ convexity and $c$ subdifferential, having
the horizontal structure of Carnot groups in mind, and we
investigate their properties. First of all, we note that in the
general situation a $c$ subdifferential is an element of the space
$\Omega_2;$ taking into account that the H--subdifferential is
contained in the first layer $V_1,$ we consider a \lq\lq
cost\rq\rq\ function
$$c:\mathbf{G}\times V_1\to\R.$$

Now, we are in the position to introduce our main definition:
\begin{definition}\label{def c H-convex}\rm
We say that a proper function $u:\Omega\subset\mathbf{G}\to
(-\infty,+\infty]$ is a \emph{$c$ H--convex function} if for every
$g\in\Omega$ we have
$$
u(g)=\sup_{(v,\alpha)\in \Pca_g}(c(g,v)+\alpha),
$$
where $ \Pca_g=\{(v,\alpha)\in V_1\times\R:\ c(g',v)+\alpha\le
u(g'),\ \forall g'\in H_{g}\cap\Omega\}$ is, for every $g\in
\Omega,$ a nonempty set.
\end{definition}
Moreover, we say that $u$ is $c$ H--concave if $-u$ is $c$
H--convex.

We would like to stress the difference between (\ref{c convexity
generale}) and Definition \ref{def c H-convex}: while, in the
former case, the index set $\Pca$ is fixed, in the latter one it
depends on the point $g.$ At first sight this difference is a
problem: as a matter of fact, in the classical case where
$\Omega_1=\Omega_2=\R^n$ and $c(x,y)=\langle x,y\rangle,$ the
pointwise supremum at every point $x\in\R^n$ of a family of affine
functions with parameters in a set $\Pca_x\subset\R^n\times\R$
depending on $x,$ can be a non convex function. However, if we
consider the case $\mathbf{G}=\H$ and
\begin{equation}\label{caso particolare c}
c(g,v)=\langle\xi_1(g),v\rangle,
\end{equation}
the notion of $c$ H--convexity corresponds to the so called \lq\lq
abstract $H$--convexity\rq\rq\ in \cite{CaPi2010} (see, in
particular, Definition 4.3); there, we proved that these functions
coincide with the H--convex ones defined in the previous section, at
least when they are real--valued. Indeed, the following holds:
\begin{proposition}\label{equivalenza: convex=abstract convex}
(see \cite{CaPi2010}, Theorem 1.1). If $u:\H\to\R$ and $c$ is as
in (\ref{caso particolare c}), then $u$ is $c$ H--convex if and
only if $u$ is H--convex.
\end{proposition}
This is one of the convictive reasons to say that our Definition
\ref{def c H-convex} is consistent.
 In the sequel, we say briefly that a function is \convex H--convex
if it is $c$ H--convex in a generic Carnot group $\G$ with cost
function $c$ as in (\ref{caso particolare c}).

Notice that, when dealing with $c$ H--convex functions, as well as
with $c$ convex functions, one has to face with the possible value
$+\infty;$ this gives rise to some difficulties when regularity
properties are required. The investigation of conditions entailing
the finiteness of a $c$ H--convex function will be the topic of
Proposition \ref{u finita in A_0} and is closely connected with
the nonemptiness of the $c$ H--subdifferential.

With further regularity on $c,$ one can hopefully find interesting
results about $c$ H--convex functions. To this purpose, in the
sequel, according to the context, some assumptions will be taken
into consideration:

\begin{itemize}

\item[$(\mathbf{c1})$] for every $p\in V_1,$ the function
$c(\cdot, p)$ belongs to $\Gamma^{1,1}_{\rm loc}(\G),$ with uniform
Lipschitz bound on $V_1;$

\item[$(\mathbf{c2})$] let $\Omega\subset\G;$ for every
$g\in\Omega$ and for all $\{v_n\}_n\subset V_1$ with
$\|v_n\|_{\frak{g}}\to +\infty,$ there exists $g'\in
H_g\cap\Omega$ such that
$$
\limsup_n(c(g',v_n)-c(g,v_n))=+\infty;
$$
\item[$(\mathbf{c3})$] let $\Omega\subset\G;$ for every
$g\in\Omega,$ the function $\X c(g,\cdot):V_1\to V_1$ is
one--to--one.
\end{itemize}
\noindent Notice that the function $c:\G\times V_1\to\R$ defined
in (\ref{caso particolare c}) fulfills all the properties above.

Let us spend a few words on the role that the above conditions on
$c$ will play in the sequel.  The regularity of $c$ expressed by
$(\mathbf{c1})$ will imply some regularity for any real--valued
$c$ H--convex function, like the local boundedness and the
horizontal differentiability almost everywhere. Condition
$(\mathbf{c2}),$ that represents a sort of horizontal
superlinearity of $c,$ will provide a link between the $c$
H--convexity of a function and the nonemptiness of its $c$
H--subdifferential at every point. Condition $(\mathbf{c3})$ will
be useful when dealing with the connection between the horizontal
derivatives of $c,$ of a $c$ H--convex function $u,$ and of its
$c$ H--subdifferential.

As in the classical setting, a concept strictly related to the $c$
H--convexity is the following:

\begin{definition}\label{cHdif}\rm
Let $u:\Omega\to (-\infty,+\infty],$ with $\Omega\subset\G.$ The
\emph{$c$ H--subdifferential} of $u$ at $g\in\Omega$ is the
(possibly empty) set
$$
\partial^c_Hu(g)=\{ p\in V_1:\; u(g')\ge u(g)+c(g',p)-c(g,p), \quad \forall g'\in
H_{g}\cap\Omega\}.
$$
\end{definition}
In particular, we say that $u$ is $c$ H--subdifferentiable at
$g_0$ if $\partial_H^cu(g_0)\neq \emptyset.$ Clearly, the \convex
H--subdifferential of a function coincides with its
H--subdifferential.

We will denote by $\partial_H^c u$ the multivalued map $g\mapsto
\partial_H^c u(g).$ To this purpose, given a multivalued map $T:\G\to \Pca(V_1),$
we recall that its  domain $\mathrm{dom}(T)$ is the set of points
$g\in\G$ for which $T(g)$ is nonempty, and the graph of $T$ is the
set  $\texttt{\rm graph}(T)=\{(g,v)\in\G\times V_1:\ g\in
\mathrm{dom}(T), v=T(g) \}.$

\begin{remark}\label{u infinita}\rm
From the definition of $\partial_H^c u$ we easily get that, if
$u(g)=+\infty,$ then $\partial_H^cu(g)\neq \emptyset$ if and only
if $u(g')=+\infty$ for every $g'\in H_g\cap \Omega.$
\end{remark}

Our next aim is to establish some results rephrasing those in
Proposition \ref{prop0} and in Theorem \ref{converse}, for the
more general case of $c$ H--convexity. As a matter of fact, under
suitable assumptions on the function $c,$  a characterization of
$c$ H--convexity via the nonemptiness of the $c$
H--subdifferential at every point can be given.

In order to prove next theorem, let us supply an extension of the
concept of H--Fenchel transform introduced in \cite{CaPi2010}. Let
$u:\Omega\to (-\infty,+\infty],$ with $\Omega\subset\G.$ The $c$
H--Fenchel transform of $u$ is the family of functions
$\{u^c_g\}_{g\in\Omega},$ where, for every $g\in\Omega,$ $
u^c_g:V_1\to [-\infty, +\infty]$ is given by
$$
u^c_g(v)=\sup_{g'\in H_g\cap\Omega}\left(c(g',v)-u(g')\right),
$$
for every $v\in V_1.$
 Notice that
\begin{equation}\label{disug}
u^c_g(v)\ge c(g',v)-u(g'),\qquad \texttt{\rm for all}\ g'\in H_g\cap
\Omega.
\end{equation}
Furthermore, $u^c_g(v)=-\infty$ for some $v\in V_1$ if and only if
$u(g')=+\infty$ for every $g'\in H_g\cap \Omega.$

The following theorem holds:
\begin{theorem}\label{converse c}
Let $u:\Omega\subset\G\to \R.$ If $\partial^c_Hu(g)\neq \emptyset$
for every $g\in \Omega,$ then $u$ is $c$ H--convex. Moreover, let us
suppose that $c$  satisfies $(\mathbf{c2})$ and
$c(g,\cdot):V_1\to\R$ is continuous, for every $g\in\Omega;$ if $u$
is $c$ H--convex, then $\partial^c_Hu(g)\neq \emptyset$ for every
$g\in \Omega.$
\end{theorem}
\textbf{Proof}: Assume that $\partial_H^cu(g)\neq \emptyset$ for
every $g\in \Omega.$  If $(v,\alpha)\in \Pca_g,$ then
\begin{equation}\label{dimo1}
u(g)\ge c(g,v)+\alpha.
\end{equation} We prove that, for every $g\in\Omega,$ the set $\Pca_g$ is
nonempty, and it contains an element $(v,\alpha)$ such that in
(\ref{dimo1}) we have an equality. Notice that $p\in
\partial^c_Hu(g)$ if and only if
\begin{equation}\label{prop51}
u(g)+u_g^c(p)=c(g,p). \end{equation} Indeed,
\begin{eqnarray*}
p\in \partial^c_Hu(g)&\Longleftrightarrow& u(g')\ge u(g)+c(g',p)-c(g,p) ,\quad\forall g'\in H_g\cap \Omega\\
&\Longleftrightarrow& c(g,p)-u(g)\ge u^c_g(p).
\end{eqnarray*}
Taking into account (\ref{disug}), we obtain that (\ref{prop51})
holds, and that $(p,-u^c_g(p))$ belongs to $\Pca_g.$ Hence, for
every $g\in\Omega,$ we have that
$u(g)=\sup_{\Pca_g}\{c(g,p)+\alpha\},$ thereby proving that $u$ is
$c$ H--convex.

Conversely, fix $g_0\in \Omega.$ Since $u$ is $c$ H--convex, there
exists a sequence $\{(p_n,\alpha_n)\}_n\subset \Pca_{g_0}$ such
that
\begin{eqnarray}
&&c(g,p_n)+\alpha_n\le u(g),\qquad \forall g\in H_{g_0}\cap
\Omega\label{dim pro3a}\\
&&c(g_0,p_n)+\alpha_n \to  u(g_0),\nonumber
\end{eqnarray}
with
\begin{equation}\label{dim pro3b}
 c(g_0,p_n)+\alpha_n -
u(g_0)>-1/n.
\end{equation}
Inequalities (\ref{dim pro3a}) and (\ref{dim pro3b}) give, for every
$g\in H_{g_0}\cap \Omega$ and for every $n,$
\begin{equation}\label{dim pro3c}
u(g)> u(g_0)+c(g,p_n)-c(g_0,p_n)-1/n.
\end{equation}
Let us first prove that $\{p_n\}_n$ is bounded in $V_1.$ By
contradiction, suppose that  $\{p_n\}_n$ is unbounded; hence, by
$(\mathbf{c2}),$ there exists $g'\in H_{g_0}\cap \Omega$ such
that, by (\ref{dim pro3c}),
$$
u(g') \ge \limsup_n (u(g_0)+c(g',p_n)-c(g_0,p_n)-1/n)=+\infty.
$$
This contradicts the assumption that $u$ is real--valued.
Therefore $\{p_n\}_n$ is bounded in $V_1$ and we can suppose that
$p_n\to p\in V_1.$ The continuity of $c(g_0,\cdot)$ and (\ref{dim
pro3b}) imply that $ \alpha_n\to -c(g_0,p)+u(g_0)$ and,
consequently,
\begin{equation}\label{dim pro3d}
\alpha_n-u(g_0)+c(g_0,p)>-1/n
\end{equation}
for sufficiently large $n.$ For every $g\in H_{g_0}\cap \Omega$ and
large $n$, (\ref{dim pro3a}), (\ref{dim pro3d}) and the continuity
of $c(g_0,\cdot)$  give
\begin{eqnarray*}
u(g)&\ge& \lim_n(c(g,p_n)+\alpha_n)\\
&\ge& \lim_n(c(g,p_n)+u(g_0)-c(g_0,p)-1/n)\\
&\ge& c(g,p)+u(g_0)-c(g_0,p)
\end{eqnarray*}
This proves that $p\in \partial_H^c u(g_0).$  \QED

A consequence of the previous result is an extension of
Proposition \ref{equivalenza: convex=abstract convex}.
\begin{corollary}\label{confronto convessità H con c}
Let $\Omega\subset\mathbf{G}$ be an H--convex, open set, and let
$u:\Omega\to\R.$ If the function $u$ is \convex H--convex, then
$u$ is H--convex. If the function $u$ is H--convex, and
measurable if $r>2,$ then $u$ is \convex H--convex.
\end{corollary}
{\bf Proof:} The proof follows from Theorem \ref{converse},
Theorem \ref{converse c} and Proposition \ref{prop0}. \QED

Next two examples show that finiteness is a binding condition for
the previous results, that fail when non real--valued functions
are involved. Consequently, the investigation about the finiteness
of a $c$ H--convex function is critical (see Proposition \ref{u
finita in A_0}).

\footnotesize
\begin{example}\label{esempio} {\rm Consider the $\R^3$--convex
function $u:\H\to (-\infty,+\infty]$ defined by
$$
u(x,y,t)=\left\{
\begin{array}
{ll} 0&t\le 0\\ +\infty &t>0 .
\end{array}
\right.
$$
It is an exercise to show that $u$ is not \convex\footnotesize
H--convex, while it is H--convex. }
\end{example}\normalsize
The previous example shows that the mentioned class of \lq\lq
abstract $H$--convex\rq\rq\ functions and the class of \convex
H--convex functions are coincident only for real--valued
functions. In next example, $c$ fulfills assumption
$(\mathbf{c2})$ together with a stronger concavity requirement,
but this seems to be irrelevant.

\footnotesize
\begin{example} {\rm Let us consider the function $u:\H\to (-\infty,+\infty]$ defined
as follows:
$$
u(x,y,t)=\cases{+\infty&\texttt{\rm if}\ $\max\{x,y\}> 0$\cr 0&
\texttt{\rm if}\ $\max\{x,y\}\le 0.$\cr}
$$
This function turns out to be $c$ H--convex, with $c(g,p)=-\langle
\xi_1(g)-p,\xi_1(g)-p\rangle;$ indeed, tedious computations show
that, for every $g\in \H,$
\begin{eqnarray*}
\Pca_g=\{((v_1,v_2),\alpha):&& \alpha\le 0\quad \texttt{\rm if}\
v_1\le 0\  \texttt{\rm and}\ v_2\le 0,\quad  \alpha\le v_1^2\
\texttt{\rm if}\ v_1> 0\  \texttt{\rm
and}\ v_2\le 0,\\
&& \alpha\le v_2^2\ \texttt{\rm if}\ v_1\le 0\  \texttt{\rm and}\
v_2>0,\quad  \alpha\le v_1^2+v_2^2\ \texttt{\rm if}\  v_1> 0\
\texttt{\rm and}\ v_2> 0\}.
\end{eqnarray*}
However, at any point $g=(x,y,t)$ such that $\max\{x,y\}>0,$  the
set $\partial_H^cu(g)$ is empty. We remark that, for every fixed
$v\in V_1,$ the function $c$ is strictly H--concave.}
\end{example}
\normalsize

In general, it is reasonable to detect some properties about $c$
implying the inclusion of the class of the H--convex functions in
the class of the $c$ H--convex functions. Next result provides a
comparison between H--convexity and $c$ H--convexity for
real--valued functions; a similar one in the classical Euclidean
case can be found in \cite{GuNg2007}, Proposition 2.4.
\begin{proposition}\label{c convex}
Assume that,  for every $p\in V_1,$ the function $c(\cdot,p)$ is
H--concave in $\G$ and, for every $g\in \Omega,$
\begin{equation}\label{con c convex}
\bigcup_{v\in V_1}\partial^H c(g,v)=V_1.
\end{equation}
Let $u:\Omega\subset\G\to\R$ be an H--convex function on the open,
H--convex set $\Omega$; moreover, if $r>2$, we assume that $u$ is
measurable. Then $u$ is $c$ H--convex. In particular, any affine
function $\phi(g)=\langle \xi_1(g),v\rangle +\alpha,$ with $v\in
V_1$ and $\alpha\in\R,$ is $c$ H--convex.
\end{proposition}
{\bf Proof:} Fix any $g\in \Omega.$ From the assumptions on $u$
and from Theorem \ref{converse}, we have that $\partial_H
u(g)\not=\emptyset.$ Hence there exists $p\in
\partial_H u(g)$ such that:
\begin{equation}\label{dim pro1}
u(g')\ge u(g)+\langle p,\xi_1(g')-\xi_1(g)\rangle, \quad\forall
g'\in H_g\cap \Omega.
\end{equation}
From (\ref{con c convex}), let $v=v(g,p)$ be such that $p\in
\partial^H c(g,v).$ By assumption we have that
\begin{equation}\label{dim pro2}
c(g',v)\le c(g,v)+\langle p,\xi_1(g')-\xi_1(g)\rangle,
\quad\forall g'\in H_g\cap \Omega.
\end{equation}
Inequalities (\ref{dim pro1}) and (\ref{dim pro2}) imply that
$v\in
\partial_H^c u(g).$ From Theorem \ref{converse c} the thesis follows. \QED

\section{Regularity properties of \emph{c} H--convex
functions}\label{section regularity}

The definition of $c$ H--convexity given in Section \ref{section
about c H-convexity}, owing to its structure,
 does not highlight any properties of
the function; in order to detect some regularity, an accurate
analysis is needed.

The problem of the regularity of a $c$ H--convex function has
already been solved for real--valued \convex H--convex functions:
indeed, Corollary \ref{confronto convessità H con c} says that a
\convex H--convex function is H--convex and hence, by the result
of Balogh and Rickly (see Theorem \ref{Ri}), it is locally
Lipschitz continuous.

In the classical situation, this investigation goes through the
notion of semiconvexity, introduced by Douglis to select unique
solutions for the Hamilton--Jacobi equation (see, for example,
\cite{Vil2009}); under suitable regularity assumptions on $c$, a
$c$ convex function $f$ is locally semiconvex, and therefore it
shares all the regularity enjoyed by convex functions (e.g., two
derivatives almost everywhere, locally Lipschitz where finite).

Encouraged by these results, we introduce the following definition:
\begin{definition}\label{locally H semiconvex} {\rm Let $\Omega\subset \mathbf{G}$ be H--convex.
A function $u:\Omega\to (-\infty,+\infty]$ is {\it H--semiconvex}
(or $\ell$ {\it H--semiconvex}) if it is $\R$--semiconvex on every
horizontal segment, i.e., there exists a positive constant $\ell$
such that
$$
u(g\delta_{\l}(g^{-1}g'))\le (1-\lambda)u(g)+\lambda
u(g')+\ell\l(1-\l)\|\xi_1(g^{-1}g')\|^2_{\frak{g}}
$$
for all $g\in\Omega$, $g'\in H_g\cap\Omega,$ and
$\lambda\in[0,1].$}
\end{definition}
We say that $u$ is locally H--semiconvex in $\Omega$ if, for every
open ball $B\subset \Omega,$ $u$ is H--semiconvex on $B;$ here and
in the sequel we consider balls arising from the gauge distance
$d$, that are H--convex. Via the equality
$$
(1-\l)\|\xi_1(g)\|^2_{\frak{g}}+\l\|\xi_1(g')\|^2_{\frak{g}}-\|\xi_1(g\delta_{\l}(g^{-1}g'))\|^2_{\frak{g}}
=(1-\l)\l\|\xi_1(g)-\xi_1(g')\|^2_{\frak{g}},
$$
an easy computation shows that $u$ is H--semiconvex if and only the
function
$$
g\mapsto u(g)+\ell \|\xi_1(g)\|^2_{\frak{g}}
$$
is H--convex. Hence, the characterization (\ref{caratterizzazione
con gradiente H}) for H--convex functions in $\Gamma^1(\Omega),$
where $\Omega$ is open, gives us that $u$ is $\ell$ H--semiconvex
if and only if
\begin{equation}\label{car Hsemiconvex}
\langle \X u(g)+2\ell\xi_1(g)-\X u(g')-2\ell\xi_1(g'),
\xi_1(g)-\xi_1(g')\rangle\ge 0,\qquad \forall g\in\Omega,\ g'\in
H_g\cap\Omega.
\end{equation}
Moreover, if $u\in \Gamma^2(\Omega),$ then $u$ is $\ell$
H--semiconvex if and only if $[\X^2 u]^*\ge -2\ell I$ within
$\Omega,$ i.e., $[\X^2 u(g)]^*+2\ell I$ is positive semidefinite
for all $g\in\Omega.$

The following fundamental proposition is the horizontal version of
a result in \cite{GaMc1996} (see Proposition C2); as a matter of
fact, our proof is completely different on account of the
definition of $c$ H--convexity:

\begin{theorem}\label{da c Hconvex a loc Hsemiconvex}
Let $ \psi:\Omega\subset\G\to(-\infty,+\infty]$ be a proper $c$
H--convex function. Assume that $(\mathbf{c1})$ is satisfied,
 i.e., for every open ball $B\subset \Omega$ there exists
$K_B>0$ such that
\begin{equation}\label{c gamma11loc}
 \|\X c(g', v)-\X c(g,v)\|_{\frak{g}}\le 2K_B d(g',g),\quad \forall g',g\in
B \ \texttt{\it and}\ \ \forall v\in V_1.
\end{equation}
Then,
 $\psi$ is locally H--semiconvex.
\end{theorem}
{\bf Proof:} Fix an open ball $B\subset \Omega.$ By the
assumptions, for all $g\in B,$ $g'\in H_g\cap B,$ and $v\in V_1,$
we have
$$
\|\X c(g', v)-\X c(g,v)\bigr\|_{\frak{g}}\le 2K_B
\|g^{-1}g'\|_{\G}=2K_B\| \xi_1(g')-\xi_1(g)\|_{\frak{g}} ,
$$
since $\xi(g^{-1}g')\in V_1.$ In particular,
$$
\langle \X c(g, v)-\X c(g',v), \xi_1(g')-\xi_1(g)\rangle\le
2K_B\langle \xi_1(g')-\xi_1(g), \xi_1(g')-\xi_1(g)\rangle.
$$
Hence, by (\ref{car Hsemiconvex}), we have that $c(\cdot,v)$ is
$K_B$ H--semiconvex in $B$ for every $v\in V_1.$

Since $\psi$ is $c$ H--convex, by definition
$$
\psi(g)=\sup_{(v,\alpha)\in \Pca_g}(c(g,v)+\alpha),
$$
where $ \Pca_g=\{(v,\alpha)\in V_1\times\R:\ c(g',v)+\alpha\le
\psi(g'),\ \forall g'\in H_{g}\cap \Omega\}. $ Let us consider the
function $\phi:B\to(-\infty,+\infty]$  defined, for every $g\in B,$
by
\begin{equation}\label{def psi}
\phi(g)=\psi(g)+K_B\|\xi_1(g)\|^2_{\frak{g}}.
 \end{equation}
 We will prove that $\phi$ is H--convex on $B.$
By contradiction, assume that there exist $g\in B$ and $g'\in
H_g\cap B$ such that $\phi$ is not $\R$--convex along the points
of the horizontal segment $\sigma_{g,g'};$ for every $\l\in[0,1],$
we denote by $g_\l$ the point $\sigma_{g,g'}(\l)$ of such
horizontal segment (see (\ref{twiset convex combination})). The
following three cases can occur:

\noindent \emph{First case:} $\phi$ is real--valued on
$\sigma_{g,g'};$ in this case, there exists $\lambda \in (0,1)$
such that
$$
2\epsilon=\phi(g_\lambda)-(\phi(g)(1-\l)+\phi(g')\l),
$$
for some positive $\epsilon.$ From the definition (\ref{def psi})
of $\phi$ and the $c$ H--convexity of $\psi,$ there exists
$(p,\alpha)\in \Pca_{g_\lambda}$ such that
$$
\phi(g_\lambda)-  (c(g_\l,p)+\alpha+K_B\|\xi_1(g_\l)\|^2_{\frak{g}})<\epsilon
$$
with
$$
c(g'',p)+\alpha\le \psi(g''),\qquad \forall g''\in
H_{g_\lambda}\cap\Omega.
$$
From
$$
c(g,p)+\alpha\le \psi(g) \quad\texttt{\rm and} \quad
c(g',p)+\alpha\le \psi(g'),
$$
and since $c(\cdot,p)+K_B\|\xi_1(\cdot)\|^2_{\frak{g}}$ is
H--convex in $B,$ we get
\begin{eqnarray*}
\phi(g_\lambda)&<& c(g_\l,p)+\alpha+K_B\|\xi_1(g_\l)\|^2_{\frak{g}}+\epsilon\\
&\le&  (1-\lambda)(c(g,p)+\alpha+K_B\|\xi_1(g)\|^2_{\frak{g}})
 +\lambda (c(g',p)+\alpha+K_B\|\xi_1(g')\|^2_{\frak{g}})+\epsilon\\
&\le&  (1-\lambda)(\psi(g)+K_B\|\xi_1(g)\|^2_{\frak{g}})+\lambda (\psi(g')+K_B\|\xi_1(g')\|^2_{\frak{g}})+\epsilon\\
&=&  (1-\lambda)\phi(g)+\lambda \phi(g')+\epsilon\\
&=& \phi(g_\lambda)-\epsilon,
\end{eqnarray*}
a contradiction.

\noindent \emph{Second case:} $\phi$ is finite at the endpoints
$g$ and $g',$ but $\phi(g_\lambda)=+\infty$ for some $\lambda\in
(0,1).$ Then, by definition of $\phi,$ for infinitely many
integers $n$ there exists $(p_n,\alpha_n)\in \Pca_{g_{\lambda}}$
such that
$$
n< c(g_\lambda,p_n)+\alpha_n+K_B\|\xi_1(g_\l)\|^2_{\frak{g}},
$$
with
$$
c(g'',p_n)+\alpha_n\le \psi(g''),\qquad \forall g''\in
H_{g_\lambda}\cap\Omega.
$$
From
$$
c(g,p_n)+\alpha_n\le \psi(g) \quad\texttt{\rm and} \quad
c(g',p_n)+\alpha_n\le \psi(g'),
$$
and since $c(\cdot,p)+K_B\|\xi_1(\cdot)\|^2_{\frak{g}}$ is
H--convex in $B,$ we get
\begin{eqnarray*}
n&<& c(g_\l,p_n)+\alpha_n+K_B\|\xi_1(g_\l)\|^2_{\frak{g}}\\
&\le&  (1-\lambda)(c(g,p_n)+\alpha_n+K_B\|\xi_1(g)\|^2_{\frak{g}})
 +\lambda (c(g',p_n)+\alpha_n+K_B\|\xi_1(g')\|^2_{\frak{g}})\\
&\le&  (1-\lambda)(\psi(g)+K_B\|\xi_1(g)\|^2_{\frak{g}})+\lambda (\psi(g')+K_B\|\xi_1(g')\|^2_{\frak{g}})\\
&=&  (1-\lambda)\phi(g)+\lambda \phi(g')\\
&\le& \max\{\phi(g),\phi(g')\},
\end{eqnarray*}
a contradiction.

\noindent \emph{ Third case:} If $\phi(g)=+\infty,$ or
$\phi(g')=+\infty,$ then (\ref{def Hconvex}) holds for every
$\lambda\in [0,1].$

\noindent Hence $\phi$ is H--convex in $B$ and the thesis follows.
\QED

In the case (\ref{caso particolare c}), condition $(\ref{c
gamma11loc})$ is satisfied with $K_B=0$ and hence, as a
consequence of the previous result, we have that
\begin{remark}\label{inclusione}{\rm
Every \convex H--convex function $u:\G\to(-\infty,+\infty]$ is
H--convex.}
\end{remark}

We note that Theorem \ref{da c Hconvex a loc Hsemiconvex} is a
very general result for proper function. In the next section, we
will investigate the problem of the finiteness of a $c$ H--convex
function. However, the proposition above and the result by Balogh
and Rickly  give rise to some interesting regularity for
real--valued $c$ H--convex functions:

\begin{corollary}\label{regolarità semiconvex} Let $\Omega$ be an open, H--convex subset of
$\G,$ and $\psi:\Omega\to\R$ a $c$ H--convex function, measurable
if $r>2.$ Assume that $(\mathbf{c1})$ holds. Then,
\begin{itemize}
 \item[i.]  $\psi$ is locally bounded;\item[ii.] $\X\psi$
exists a.e. on every open ball $B\subset\Omega .$
\end{itemize}
\end{corollary}
{\bf Proof:} \emph{i.} From Theorem \ref{da c Hconvex a loc
Hsemiconvex}, for every open ball $B\subset\Omega$ there exists
$\ell>0$ such that the function $g\mapsto
\Psi(g)=\psi(g)+\ell\|\xi_1(g)\|^2_{\frak{g}}$ is H--convex on
$B,$ and, by the assumptions, it is measurable if $r>2.$ From
Theorem \ref{Ri}, since $\Psi$ is Lipschitz on every ball
$B\subset\Omega,$ $\Psi$ is bounded on $B;$ this implies the
boundedness of $\psi$ on $B,$ for every $B\subset \Omega.$

\noindent  \emph{ii.}  From Theorems \ref{Ri} and \ref{DaGaSa},
$\X\Psi(g)$ exists for almost all $g\in B;$ we conclude that
$\X\psi(g)$ exists for almost all $g\in B.$\QED

In the Euclidean case (see, for instance, Proposition 2.7 in
\cite{GuNg2007}) a connection can be stated between the $c$
subdifferential of a function $f,$ and the gradients $\nabla c$ and
$\nabla f;$ a perfectly symmetrical result holds in our framework.

\begin{proposition}\label{partial_H^c u e partial_H u}
Let $c:\mathbf{G}\times V_1\to \R$ and $u:\Omega\subset \G\to\R$ be
such that $\X u(g_0)$ and $\X c(g_0,v)$ exist for every $v\in V_1$
and for some $g_0\in \mathrm{int}(\Omega).$
\begin{itemize}
\item[i.] If $p\in
\partial_H^cu(g_0),$ then
$ p\in\left( \X c(g_0,\cdot)\right)^{-1}(\X u(g_0)); $

\item[ii.] if $\partial_H^c u(g_0)\neq \emptyset$ and $\X
c(g_0,\cdot):V_1\to V_1$ is one--to--one, then
\begin{equation}\label{p esplicito}
\partial_H^cu(g_0)=\left\{\left( \X c(g_0,\cdot)\right)^{-1}(\X u(g_0))\right\}.
\end{equation}
\end{itemize}
\end{proposition}

\noindent{\bf Proof:} From the definition of $c$
H--subdifferential, for all $g\in H_{g_0}\cap \Omega,$ we have
that
$$
u(g)-c(g,p)\ge u(g_0)-c(g_0,p);
$$
in particular, $g_0$ is a minimum point for the function $g\mapsto
u(g)-c(g,p)$ on the plane $H_{g_0}\cap\Omega.$ This implies that
\begin{equation}\label{uguaglianza gradienti}
\X  u(g_0)=\X  c(g_0, p).
\end{equation}
If we consider the function $\X  c(g_0,\cdot):V_1\to V_1,$ then
(\ref{uguaglianza gradienti}) implies \emph{i}. The additional
assumption in \emph{ii}. gives easily (\ref{p esplicito}). \QED

Under more regularity assumptions on $c$, Theorem \ref{converse
c}, Corollary \ref{regolarità semiconvex} and Proposition
\ref{partial_H^c u e partial_H u} entail the following
\begin{corollary}\label{corollario fine section 5}
Let $\Omega$ be an open, H--convex subset of $\G,$ and
$u:\Omega\to\R$ be a $c$ H--convex function, measurable if $r>2.$
Assume that $c$ fulfills $(\mathbf{c1}),$ $(\mathbf{c2})$ and
$(\mathbf{c3}),$ and that $c(g,\cdot):V_1\to\R$ is continuous, for
every $g\in\Omega.$

Then, for a.e. $g\in \Omega,\ \partial_H^c u(g)$ is a singleton,
and $\partial_H^c u(g)=\{(\X c(g,\cdot))^{-1}(\X u(g))\}.$
\end{corollary}

Exploiting the previous results, a necessary condition for a
function to be $c$ H--convex can be given. Let $\Omega$ be an open
set, and $c$ be a function satisfying the assumptions of Corollary
\ref{corollario fine section 5}; we assume, in addition, that
$c(\cdot, v)\in \Gamma^2(\Omega),$ for every $v\in V_1.$ Consider
a $c$ H--convex function $u\in \Gamma^2(\Omega);$ then, from
Theorem \ref{converse c} and Proposition \ref{partial_H^c u e
partial_H u}, we get that $\partial_H^c u(g)$ is a singleton and
it is given by (\ref{p esplicito}). For every $g_0\in\Omega,$
denote by $p_0$ the unique $c$ H--subgradient of $u$ at $g_0;$
then, the function
$$
\theta_{g_0}:H_{g_0}\cap\Omega\to\R,\quad
\theta_{g_0}(g)=u(g)-c(g,p_0),
$$
has a minimum at $g_0.$ This implies that $ [\X^2
\theta_{g_0}]^*(g_0)\ge 0.$
From (\ref{p esplicito}), we obtain a
necessary condition for the $c$ H--convexity of $u:$
$$
[\X^2 u]^*(g)\ge  [\X^2 c]^*\biggl(g,(\X c(g,\cdot))^{-1}(\X
u(g))\biggr), \quad \forall g\in \Omega.
$$
In the particular
situation where $c(g,v)=-\|\xi_1(g)-v\|^2_{\frak{g}},$ we obtain
$$
[\X^2 u]^*(g)\ge -2I, \quad\forall g\in \Omega.
$$

\section{\emph{c} H--cyclic monotonicity}\label{$c$ H--cyclic monotonicity}
In $\R^n$ and, more generally, in Banach spaces $\mathbf{X}$,
the graph of the multivalued map defined via the subdifferential
$\partial f$ of a function $f$ is a cyclically monotone subset of
$\mathbf{X}\times \mathbf{X}^*,$ i.e.
$$
\sum_{i=0}^n \langle x_{i+1},x^*_i\rangle \le \sum_{i=0}^n \langle
x_{i},x^*_i\rangle,
$$ for every finite sequence $\{(x_i,x_i^*)\}_{i=0}^n\subset \texttt{\rm
graph}(\partial f),$ with $x_{n+1}=x_0.$  A cyclically monotone
subset in $\mathbf{X}\times \mathbf{X}^*$ is called maximal if it
is not a proper subset of another cyclically monotone set in
$\mathbf{X}\times \mathbf{X}^*.$ In this context, a well--known
result due to R.T. Rockafellar \cite{Ro1970} says that the maximal
cyclically monotone subsets of $\mathbf{X}\times \mathbf{X}^*$ are
completely characterized as the graphs of the multivalued maps
$x\mapsto
\partial f(x),$ where $f$ is a proper lower semicontinuous convex
function.

This result was extended to the case of $c$ convex functions
$f:\Omega_1\to(-\infty,+\infty]$ and $c$ cyclically monotone sets
$\Gamma\subset \Omega_1\times \Omega_2,$ where $\Omega_i$ are very
general spaces  (see, for instance, \cite{Ru1996}). We recall that
$\Gamma$ is said to be $c$ cyclically monotone if for all $\{(x_i,
y_i)\}_{i=0}^n\subset \Gamma,$ with $x_{n+1}=x_0,$
\begin{equation}\label{c cyclicall classica}
\sum_{i=0}^nc(x_{i+1},y_i)\le \sum_{i=0}^nc(x_{i},y_i).
\end{equation}
We would like to stress that the $c$ subdifferential of $f$ at a
point is a (possibly empty) subset of $\Omega_2$ defined in
(\ref{c subdifferential generale}).

The aim of this section is to adapt Rockafellar's ideas in
\cite{Ro1966} to the sub--Riemannian structure of a Carnot group,
in the \lq\lq$c$ case\rq\rq. Two are the main features of our
setting. First, the horizontal subdifferential $\partial_Hu,$ that
plays a fundamental role in the study of the horizontal convexity
of $u,$ is a subset of $V_1.$ Thereby the graph of the map
$g\mapsto
\partial_Hu(g)$ is a subset of $\G\times V_1,$ and this is the main
reason why we will introduce the notion of $c$ H--cyclic
monotonicity for a subset of $\mathbf{G}\times V_1.$ Furthermore,
the H--subdifferential of a function at a point carries
information about the function only along horizontal segment
through the point itself. To this purpose, in \cite{CaPi2008} we
proved that, if $u$ is a real--valued, H--convex function on $\H,$
and so H--subdifferentiable, then their H--subgradients are
sufficient to \lq\lq reconstruct\rq\rq\ the function. More
precisely, using the definition of H--sequence (see Section
\ref{section H-convexity}), we proved the following
\begin{theorem}\label{Teo CaPi2010} (see \cite{CaPi2008}, Theorem 6.4). If  $u:\H\to\R$ is an H--convex function, then
\begin{equation}\label{Rock function}
u(g)=u(g_0)+\sup_{{\cal Q}_g} \left\{\sum_{i=0}^{n-1} \langle
p_i,\xi_1(g_{i+1})-\xi_1(g_i)\rangle\right\},
\end{equation}
where $g_0$ is fixed, and ${\cal
Q}_g=\bigl\{\{(g_i,p_i)\}_{i=0}^{n}: \{g_i\}_{i=0}^{n}\
\texttt{\it H--sequence},\ g_n=g,\ p_i\in\partial_Hu(g_i) \bigr\}.
$
\end{theorem}
In the sequel, we deal with the more general case of $c$ H--convex
functions. To begin, let us investigate about the finiteness of a
$c$ H--convex function.

First of all, given a subset  $A$ of $\G$ and a point $g_0\in A,$
we will consider a particular set of points that has a good
behaviour with respect to horizontal displacements from $g_0$
within $A.$ Let us denote by $\mathcal{H}(g_0,A)$ the subset of
$A$ that contains exactly those points that can be reached
starting from $g_0$ and moving along horizontal segments whose
endpoints lye in $A.$ More precisely, a point $g$ belongs to
$\mathcal{H}(g_0,A)$ if there exists an H--sequence $\{
g_i\}_{i=0}^{n}$ such that $g_n=g$ and $g_i\in A,$ for every
$i=1,2,\dots, n.$ In some cases, this set is a singleton; as an
example, if $A=\{(0,0,t)\in \H: t\in\R\}$ and $g_0=(0,0,0),$ we
get $\mathcal{H}(g_0,A)=\{g_0\}.$ From Proposition \ref{FoSt} we
easily get the following
\begin{remark}\label{interior point}\rm
If $g_0$ is an interior point of $A$, then $g_0$ is an interior
point of $\mathcal{H}(g_0,A).$
\end{remark}

In the next proposition, we prove a sufficient condition for the
finiteness of a $c$ H--convex function defined on a set $\Omega$,
at least on $\mathcal{H}(g_0,\Omega).$ This result will play a
fundamental role in the main theorem of this section. Let us
recall that, for a given multivalued map $T,$
$\mathrm{dom}(T)=\{g\in \G:\,T(g)\neq \emptyset\}.$

\begin{proposition}\label{u finita in A_0}
Let $u:\Omega\subset\G\to(-\infty,+\infty],$ and let $g_0\in
\Omega$ be such that $u(g_0)<+\infty$ and $\partial_H^c u(g_0)\neq
\emptyset.$ Then, $u$ is real--valued in
${\mathcal{H}(g_0,\mathrm{dom}(\partial_H^c u))}.$
\end{proposition}
{\bf Proof:} For any $g\in
\mathcal{H}(g_0,\mathrm{dom}(\partial_H^c u)),$  $g\neq g_0,$
there exists an H--sequence $\{g_i\}_{i=0}^{n}$ with $g_n=g$ and
$\partial_H^c u(g_i)\neq \emptyset,$ for every $i=0,1,\dots, n.$
This implies that $u(g)\neq +\infty.$ Indeed, since $g_i$ and
$g_{i+1}$ are the endpoints of a horizontal segment, from Remark
\ref{u infinita} it follows that $u(g_i)=+\infty$ if and only if
$u(g_{i+1})=+\infty.$ Since $u(g_0)\neq +\infty,$ we get the
result. \QED

Recalling that an H--sequence $\{g_i\}_{i=0}^{n}$ is closed when
$g_n\in H_{g_0}$ (in this case we set $g_{n+1}=g_0$), we give the
following natural
\begin{definition}\label{definizione c Hcyclically monotone}\rm
We say that $\mathcal{R}\subset \mathbf{G}\times V_1$ is a $c$
H\emph{--cyclically monotone set} if, for every sequence
$\{(g_i,p_i)\}_{i=0}^{n}\subset \mathcal{R}$ such that
$\{g_i\}_{i=0}^{n}$ is a closed H--sequence, we have that
\begin{equation}\label{def c Hcyclically}
\sum_{i=0}^{n} c(g_{i+1},p_i)\le \sum_{i=0}^{n} c(g_{i},p_i).
\end{equation}
\noindent We say that a multivalued map $T:\mathbf{G}\to
\mathcal{P}( V_1)$ is a $c$ H\emph{--cyclically monotone map} if
$\texttt{\rm graph}(T)$ is $c$ H--cyclically monotone.
\end{definition}

From this definition, we can express in a different way the
characterization of H--convex functions in $\Gamma^1$ presented in
(\ref{caratterizzazione con monotomia gradiente H}):
\begin{remark}\rm
Let $u\in\Gamma^1(\Omega).$ Then $u$ is H--convex if and only if the
map $g\mapsto
\partial_H u(g)$ has a \convex H--cyclically monotone graph.
\end{remark}
Notice that, for every functions $u$ and $c,$ without any
regularity assumptions, the map $g\mapsto
\partial_H^c u(g)$ has a $c$ H--cyclically monotone graph. Indeed,
if $\{(g_i,p_i)\}_{i=0}^{n}\subset\texttt{\rm graph}(\partial_H^c
u)$ and $\{ g_i\}_{i=0}^{n}$ is a closed H--sequence, then
$$
u(g_{i+1})-u(g_i)\ge c(g_{i+1},p_i)-c(g_i,p_i),\qquad i=0, \ldots,
n,
$$
implies (\ref{def c Hcyclically}).

The following result is the converse of the previous note, and it
provides a crucial link between our approach and some possible
application in optimal mass transportation problems:
\begin{theorem}\label{theorem rockafellar c}
Let $T:\G\to \mathcal{P}(V_1)$ be a $c$ H--cyclically monotone
map. Then, for all $g_0\in \mathrm{int}(\mathrm{dom}(T)),$ there
exists a $c$ H--convex function
$f_{g_0}:\mathcal{H}(g_0,\mathrm{dom}(T))\to \R, $ such that
\begin{equation}\label{inclusione teo rock}
T(g)\subset \partial_H^cf_{g_0}(g),\qquad\texttt{\it for every }\
g\in \mathcal{H}(g_0,\mathrm{dom}(T)). \end{equation}
\end{theorem}
Let us first make some comments. The function $f_{g_0},$ that will
be defined in (\ref{def funz}), is the $c$ version of the
Rockafellar's function (\ref{Rock function}) in the
sub--Riemannian setting.

The reader will infer that the function $f_{g_0},$ with $g_0\in
\mathrm{int}(\mathrm{dom}(T)),$ could be defined, using (\ref{def
funz}), at every point $g$ linked via a horizontal segment to a
point in $\mathcal{H}(g_0,\mathrm{dom}(T))$; however, one cannot
guarantee that $f_{g_0}$ is real--valued at $g$ and, above all,
that it is $c$ H--convex. Moreover, if $g_0\in \mathrm{dom}(T)$
and there does not exist any point $g\in H_{g_0}\cap
\mathrm{dom}(T),$ different from $g_0,$ using (\ref{def funz}), we
obtain a trivial function whose domain is $\{g_0\}$ and
$f_{g_0}(g_0)=-\infty.$

The domain $\mathcal{H}(g_0,\mathrm{dom}(T))$ of the function
$f_{g_0}$ would have a very strange shape. However, from Remark
\ref{interior point}, if $g_0$ is an interior point of
$\mathrm{dom}(T),$ then $g_0$ is an interior point of
$\mathrm{dom}(f_{g_0}).$

In the sequel, we will denote by $\mathrm{dom}(f_{g_0})$ the set
$\mathcal{H}(g_0,\mathrm{dom}(T)).$

\noindent{\bf Proof of Theorem \ref{theorem rockafellar c}:} Let
us suppose that $T$ is $c$ H--cyclically monotone, and fix $g_0\in
\mathrm{int}(\mathrm{dom}(T)).$ For every $g\in
\mathcal{H}(g_0,\mathrm{dom}(T))$ we define $\mathcal{Q}_g$ as the
set of all sequences $\{(g_i,p_i)\}_{i=0}^n,$ where
$\{g_i\}_{i=0}^n$ is an H--sequence with starting point $g_0,$
$p_i\in T(g_i)$ for $0\le i\le n,$ and $g_n\in H_g.$ Let
$f_{g_0}:\mathcal{H}(g_0,\mathrm{dom}(T))\to(-\infty,+\infty]$ be
the function defined by
\begin{equation}\label{def funz}
f_{g_0}(g)=\sup_{\mathcal{Q}_g}\left\{ \sum_{i=0}^{n-1}
(c(g_{i+1},p_i)-c(g_{i},p_i))+ c(g,p_n)-c(g_{n},p_n)\right\}.
\end{equation}
First of all, since $g_0\in \mathrm{int}(\mathrm{dom}(T)),$ the
set $\mathcal{Q}_g$ is nonempty, and then $f_{g_0}(g)$ is greater
than $-\infty.$ Let us show that $f_{g_0}$ is proper. For every
$g_1\in H_{g_0}\cap \mathcal{H}(g_0,\mathrm{dom}(T)),$ we have
$$
f_{g_0}(g_0)\ge c(g_1,p_0)-c(g_0,p_0)+c(g_0,p_1)-c(g_1,p_1);
$$
if we choose $g_1=g_0,$ then we obtain $f_{g_0}(g_0)\ge 0.$ Since
$T$ is $c$  H--cyclically monotone, we have that
$$\sum_{i=0}^{n-1} (c(g_{i+1},p_i)-c(g_{i},p_i))+ c(g_0,p_n)-c(g_{n},p_n)\le 0$$
for every sequence in $\mathcal{Q}_{g_0}:$ clearly this implies
$f_{g_0}(g_0)\le 0.$ Hence $f_{g_0}(g_0)=0$ and $f_{g_0}$ is
proper.

Next, let us choose $\overline{g}\in
\mathcal{H}(g_0,{\mathrm{dom}}(T)),$ and $\overline{p}\in
T(\overline{g}).$ For every
 $\alpha<f(\overline{g}),$
there exists a sequence $\{(g_i,p_i)\}_{i=0}^{n}$ in
$\mathcal{Q}_{\overline{g}}$ such that
$$
\alpha<\sum_{i=0}^{n-1} (c(g_{i+1},p_i)-c(g_{i},p_i))+
c(\overline{g},p_n)-c(g_{n},p_n).
$$
Let $g\in H_{\overline{g}}\cap \mathcal{H}(g_0,\mathrm{dom}(T)).$
By adding to the sequence above the point $(\overline{g},
\overline{p}),$ we obtain a new sequence that belongs to
$\mathcal{Q}_{{g}}.$ Then, by (\ref{def funz}),
 we have
\begin{eqnarray*}
f_{g_0}(g)&\ge& \sum_{i=0}^{n-1} (c(g_{i+1},p_i)-c(g_{i},p_i))+
c(\overline{g},p_n)-c(g_{n},p_n)+
c(g,\overline{p})-c(\overline{g},\overline{p}) \\
&>&\alpha+c(g,\overline{p})-c(\overline{g},\overline{p}).
\end{eqnarray*}
 Since $\alpha<f_{g_0}(\overline{g})$ is
arbitrary, we conclude that $\overline{p}\in
\partial_H^c f_{g_0}(\overline{g}).$ Hence we obtain (\ref{inclusione teo
rock}).

From Proposition \ref{u finita in A_0}, since $\partial_H^c
f_{g_0}({g})\not=\emptyset$ for every ${g}\in
\mathcal{H}(g_0,\mathrm{dom}(T)),$ and $f_{g_0}(g_0)$ is finite,
we can conclude that $f_{g_0}$ is real--valued in
$\mathcal{H}(g_0,\mathrm{dom}(T)).$ Finally, Theorem \ref{converse
c} and the nonemptiness of $\partial_H^c f_{g_0}(g)$ for every
$g\in \mathcal{H}(g_0,\mathrm{dom}(T))$ implies that $f_{g_0}$ is
$c$ H--convex.
 \QED

With some regularity assumptions on the function $c,$ a $c$
H--cyclic monotone multivalued map is, in fact, an a.e.
single--valued map in its domain; furthermore, the graph of $T$
coincides, locally, with the graph of the $c$ H--subdifferential
of a real--valued $c$ H--convex function. At first sight, this
seems to be a local conclusion in $\mathrm{dom}(T),$ but the
different functions $f_g$ that we construct on the sets
$\mathcal{H}(g,\mathrm{dom}(T))$, with $g\in
\mathrm{int}(\mathrm{dom}(T)),$ share indeed the same $c$
H--subdifferential. This is the content of the following
proposition that provides the sub--Riemmanian version of the
results in \cite{GaMc1996}.

\begin{corollary}\label{corollario fine section 6}
Assume that $c$ satisfies $(\mathbf{c1}),$ $(\mathbf{c2})$ and
$(\mathbf{c3}),$ and that $c(g,\cdot):V_1\to\R$ is continuous, for
every $g\in\G.$ Let $T:\G\to \mathcal{P}(V_1)$ be a $c$
H--cyclically monotone map, and denote by $g_0$ an interior point
of $\mathrm{dom}(T).$

Then there exists a real--valued $c$ H--convex function $f_{g_0}$
with the following properties:
\begin{itemize}
\item[i.] $g_0$ is an interior point of $\texttt{\rm
dom}(f_{g_0});$ \item[ii.] for every $g\in \texttt{\rm dom}(\X
f_{g_0}),$  $T(g)=\partial_H^cf_{g_0}(g)=\{\X c(g,\cdot)^{-1}(\X
f_{g_0}(g))\};$ \item[iii.] $ \texttt{\rm dom}(f_{g_0})\setminus
\texttt{\rm dom}(\X f_{g_0})$ has null measure, with the
additional assumption that $f_{g_0}$ is measurable if $r>2.$
\end{itemize}
Let $g_1$ be another point in the interior of $\mathrm{dom}(T).$
Then
\begin{itemize}
\item[iv.]
if $g\in \texttt{\rm dom}( f_{g_0})\cap\texttt{\rm dom}(
 f_{g_1}),$
we have that $T(g)\subset \partial_H^c f_{g_0}(g)\cap\partial_H^c
f_{g_1}(g);$
 \item[v.] if $g\in \texttt{\rm dom}(\X f_{g_0})\cap\texttt{\rm dom}(\X
 f_{g_1}),$
we have that
$$
\partial_H^cf_{g_0}(g)=\{\X c(g,\cdot)^{-1}(\X
f_{g_0}(g))\}=\{\X c(g,\cdot)^{-1}(\X
f_{g_1}(g))\}=\partial_H^cf_{g_1}(g).
$$
\end{itemize}
\end{corollary}

\noindent{\bf Proof:} Clearly $f_{g_0}$ is defined in Theorem
\ref{theorem rockafellar c} and consequently if finite and $c$
H--convex. Remark \ref{interior point} guarantee that \emph{i.}
holds. From the construction of the function $f_{g_0}$ in the proof
of Theorem \ref{theorem rockafellar c}, we have that $\partial_H^c
f_{g_0}(g)\not=\emptyset$ for every $g\in \texttt{\rm
dom}(f_{g_0});$ this argument and Proposition \ref{partial_H^c u e
partial_H u} imply \emph{ii}. Corollary \ref{corollario fine section
5} implies \emph{iii}. The last part of the Corollary follows from
the previous implications and Theorem \ref{theorem rockafellar c}.
 \QED

\section{An elementary application to optimal mass transportation
in $\H$}\label{application}

Recently, as we mentioned, some papers have been devoted to the
study of optimal mass transportation within Carnot groups. Whereas
it should be clear to the reader that the focus of this paper is
not this one, we would like to show, following timidly the line of
the paper by Gangbo and McCann \cite{GaMc1996}, how the tools
introduced in the previous sections could be applied, at least if
$\G=\H.$

Let $(\Omega_1,\mu)$ and $(\Omega_2,\nu)$ be probability spaces,
and let us denote by $\Gamma(\mu,\nu)$ the set of the probability
measures $\gamma$ on $\Omega_1\times \Omega_2$ with marginals
$\mu$ and $\nu$, i.e. such that $ \gamma(A\times \Omega_2)=\mu(A)$
and $\gamma(\Omega_1\times B)=\nu(B),$ for all $\mu$--measurable
sets $A$ and $\nu$--measurable sets $B.$ We say that a map
$s:\Omega_1\to \Omega_2$ pushes $\mu$ forward to $\nu,$ i.e.,
$\nu=s_\sharp\mu,$ if $\nu(B)=\mu(s^{-1}(B))$ for all
$\nu$--measurable sets $B.$

Monge's problem, formulated in 1781, takes into consideration
$\Omega_1=\Omega_2=\R^n,$ two measures $\mu$ and $\nu$ on $\R^n,$
a cost function $c:\R^n\times\R^n\to\R,$ and
\begin{equation}\label{monge}
 \inf_{\{s:\ s_\sharp
\mu=\nu\}}\int_{\R^n}c(x,s(x))\,d x.
\end{equation}
A function $s_*:\R^n\to\R^n,$ which minimizes (\ref{monge}), is
called optimal map. In 1942, Kantorovich provided a relaxed
version of the previous problem, as follows:
\begin{equation}\label{kantorovich}
 \inf_{\ \gamma\in
\Gamma(\mu,\nu)}\int_{\R^n\times\R^n}c(x,y)\,d\gamma( x,y).
\end{equation}
A measure $\gamma_*\in\Gamma(\mu,\nu) $ is an optimal measure if
it is a minimum in (\ref{kantorovich}). Since, for every $\mu$
such that $s_{\sharp}\mu=\nu,$ the measure
$\gamma=(\mathbf{1}\times s)_\sharp\mu$ belongs to $\Gamma(\mu,
\nu),$ the change of variable shows that the functional in
(\ref{kantorovich}) coincides with the one in (\ref{monge}); this
implies that the Kantorovich's infimum encompasses a large class
of objects than that of Monge.

Among the other results,  Gangbo and McCann proved that for a cost
$c(x,y)=h(x-y),$ where $h$ is a strictly convex and superlinear
function (here, for simplicity, we assume $h\in C^1$), satisfying
a technical condition that they call (H2) (see \cite{GaMc1996}, p.
121), there exists a unique solution for both the Monge and the
Kantorovich problems. In particular, if $\mu$ and $\nu$ are Borel
measures on $\R^n$ such that $\mu$ is absolutely continuous with
respect to the Lebesgue measure, and if the infimum in
(\ref{kantorovich}) is finite, then there exists a unique optimal
measure $\gamma_*=(\mathbf{1}\times s_*)_\sharp \mu,$ where $s_*$
is an optimal map for the Monge's problem that is $\mu$--a.e.
defined through a $c$ concave function $\varphi,$ usually called
\lq\lq potential\rq\rq, via the formula $s_*(x)=x-(\nabla
h)^{-1}(\nabla\varphi(x))$ (see Theorems 1.2 and 3.7 in
\cite{GaMc1996}). Here, $c$ convexity, and hence $c$ concavity,
are defined as in (\ref{c convexity generale}).

The main ingredients of this result can be summarized as follows:
if $\gamma_*$ is optimal, then its support supp($\gamma_*$) is
$(-c)$ cyclically monotone (according to (\ref{c cyclicall
classica}), with the obvious changes of the sign due to the $c$
concavity of $\varphi$). Consequently, there exists a $c$ concave
and Rockafellar's function $\varphi$ such that
$\mathrm{supp}(\gamma_*)\subset\texttt{\rm
graph}(\partial^c\varphi).$ Since $\varphi$ is locally
semiconcave, it is differentiable a.e. where it is finite; in
particular, if $x\in \texttt{\rm dom}(\nabla \varphi),$ then the
$c$ superdifferential is a singleton and it is given by
$\{x-(\nabla h)^{-1}(\nabla \varphi(x))\}.$ Finally, the function
$s_*$ defined a.e. by the condition $(x,s_*(x))\in\texttt{\rm
graph}(\partial^c\varphi),$ provides the optimal map. One moment's
reflection shows that the mentioned objects and tools have already
been defined in the previous sections in our framework.

Let $(\H,\mu)$ and $(\H,\nu)$ be probability spaces; given a
function $c:\H\times V_1\to\R,$ we define the \lq\lq profit\rq\rq\
function $C:\H\times\H\to [-\infty,\infty)$ as follows:
$$
C(g,g')= \left\{
\begin{array}{ll}
c(g,\xi_1(g'))&\texttt{\rm if}\ (g,g')\in \Sca\\
-\infty& \texttt{\rm if}\ (g,g')\not\in \Sca,
\end{array}
\right.
$$
where  $\Sca$ denotes the (symmetric) set
$\Sca=\{(g,g')\in\H\times\H:\ g'\in H_g\}.$ We study the problem
\begin{equation}\label{probsup}
 \sup_{\gamma\in\Gamma(\mu,\nu)} \Cca(\gamma),\qquad \texttt{\rm
 where}\quad
\displaystyle\Cca(\gamma)=\int_{\H\times\H} C(g,g')\, d\gamma(g,g').
\end{equation}
We say that $\gamma_*$ is optimal if
$\Cca(\gamma_*)\ge\Cca(\gamma),$ for every
$\gamma\in\Gamma(\mu,\nu).$ It is noteworthy that, with this type
of profit function, any optimal map $s_*$ moves every points, at
least a.e., along their horizontal planes,  i.e.
$s_*(g)=g\exp{v},$ for some $v=v(g)\in V_1.$ We will denote by
$S_{\gamma}$ the set $(\mathbf{1}\times \xi_1)({\rm
supp}(\gamma))\subset \H\times V_1.$

The aim of this section is to show that, for our elementary
problem (\ref{probsup}), it can be reasonably introduced a notion
of \lq\lq potential\rq\rq\ on $\H$ that identifies the optimal
map. In order to do this, we have the following:
\begin{proposition}\label{main result}
Let  $c:\H\times V_1\to (-\infty, 0]$  be a continuous function.
Let $\gamma_*$ be an optimal solution for problem (\ref{probsup}),
with $\Cca(\gamma_*)>-\infty,$ and suppose that
$\gamma'_*=(\mathbf{1}\times \xi_1)_{\#}\gamma_*$ is optimal for
\begin{equation}\label{probsup'}
 \sup_{\gamma'\in\Gamma(\mu,\nu')}\int_{\H\times V_1} c(g,v)\,
d\gamma'(g,v),
\end{equation}
where $\nu'=(\xi_1)_\sharp \nu.$ Then, the set $S_{\gamma_*} $ is
$c$ H--cyclically monotone.
\end{proposition}
The assumptions of the proposition above deserve some comments.
Indeed, for any $\gamma\in \Gamma(\mu,\nu),$ the measure
$\gamma'=(\mathbf{1}\times \xi_1)_\sharp\gamma$ is in
$\Gamma(\mu,\nu')$. On the contrary, if $\gamma'$ belongs to
$\Gamma(\mu,\nu'),$ one cannot infer, in general, the existence of
$\gamma\in \Gamma(\mu,\nu)$ such that $(\mathbf{1}\times
\xi_1)_\sharp\gamma=\gamma'.$ This implies that, if $\gamma_*$ is
optimal for (\ref{probsup}), one cannot infer that
$(\mathbf{1}\times \xi_1)_\sharp\gamma_*$ is optimal for
(\ref{probsup'}). \psn

\noindent{\bf Sketch of the proof:} First of all notice that, if
$c$ is bounded from above, then for any $\gamma\in\Gamma(\mu,\nu)$
such that $\Cca(\gamma)>-\infty ,$  supp$(\gamma)\subset \Sca.$ By
the change of variables theorem, we get
$$
\Cca(\gamma)=\int_{\Sca}
C(g,g')\, d\gamma(g,g')\\
=\int_{S_\gamma} c(g,v)\, d\gamma'(g,v),
$$
where $\gamma'=(\mathbf{1}\times \xi_1)_\sharp\gamma.$

Let $\gamma_*$ and $\gamma'_*$ satisfy the assumptions. By
contradiction, assume that $S_{\gamma_*}$ is not $c$ H--cyclical
monotone; then, there exists $\{(g_i^*,p_i^*)\}_0^n\subset
\mathrm{supp}(\gamma'_*),$ where $\{g_i^*\}_0^n$ is a closed
$H$--sequence, such that the continuous function $f:\H^{n+1}\times
V_1^{n+1}\to \R$
$$
f(g_0,g_1,\dots, g_n,p_0,p_1,\dots,
p_n)=\sum_{i=0}^n\left(c(g_{i+1},p_i)-c(g_i,p_i)\right)
$$
is positive at $g_i=g_i^*$ and $p_i=p_i^*.$ At this step, the
proof follows the same line of Theorem 2.3 in \cite{GaMc1996},
showing that $\gamma'_*$ cannot be an optimal measure for problem
(\ref{probsup'}). \QED

The result above allows us to connect the tools of the previous
sections to the optimal transportation, and to introduce a notion
of \lq\lq potential\rq\rq\ in the Heisenberg framework. Since this
will be defined via the Rockafellar's function of Theorem
\ref{theorem rockafellar c}, we must take into account that such
theorem provides only local information.

Let $c:\H\times V_1\to (-\infty, 0]$ be a continuous function
satisfying $(\mathbf{c1}),$ $(\mathbf{c2})$ and $(\mathbf{c3}),$
and $\gamma_*$ and $\gamma'_*$ be as in Proposition \ref{main
result}. We consider the multivalued map $T_{\gamma_*}:\H\to
\mathcal{P}(V_1)$ defined as $$T_{\gamma_*}(g)=\{v\in V_1:\
(g,v)\in S_{\gamma_*}\}.
$$
Proposition \ref{main result} guarantees that $T_{\gamma_*}$ is a
$c$ H--cyclically monotone map. From Corollary \ref{corollario
fine section 6}, there exists a family of $c$ H--convex functions
$$\Psi=\{\psi_g:\ g\in\texttt{\rm int}(\texttt{\rm dom}(T_{\gamma_*}))\}$$
such that $\texttt{\rm dom}(\psi_g)$  is a subset of $\texttt{\rm
dom}(T_{\gamma_*})$ and contains $g$ as an interior point.
Moreover, for every $g\in \texttt{\rm int}(\texttt{\rm
dom}(T_{\gamma_*}))$ and for a.e. $g'\in \texttt{\rm dom}(
\psi_g),$ there exists $\X \psi_g(g')$ and hence $\partial_H^c
\psi_g(g')$ is a singleton. Finally, if $g'$ is in the domain of
two functions $\psi_{g_1},\ \psi_{g_2}$ in $\Psi,$ then
$\partial_H^c \psi_{g_1}(g')$ and $\partial_H^c \psi_{g_2}(g')$
have nonempty intersection.

For these reasons, given a point $g\in \texttt{\rm
int}(\texttt{\rm dom}(T_{\gamma_*})),$ we define the $c$
H--subdifferential of the family $\Psi$ at $g$ as the set
$$
\partial_H^c \Psi(g)=\bigcap_{\{g'\in {\rm int}({\rm dom}(T_{\gamma_*})):\ g\in {\rm dom}(\psi_{g'})\} }
\partial_H^c \psi_{g'}(g).$$
Since, for every $g'\in \texttt{\rm int}(\texttt{\rm
dom}(T_{\gamma_*}))$ and for every $g\in  {\rm dom}(\psi_{g'}),$
Theorem \ref{theorem rockafellar c} guarantees that
$T_{\gamma_*}(g)\subset
\partial_H^c \psi_{g'}(g),$ we have that $
\partial_H^c \Psi(g)$ is nonempty. Clearly, for a.e. $g\in\texttt{\rm
int}(\texttt{\rm dom}(T_{\gamma_*})),$ the set $\partial_H^c
\Psi(g)$ is a singleton and it defines a.e. the optimal map $s_*.$
More precisely, if $\partial_H^c \Psi(g)$ is a singleton, then
$(g,\xi_1(s_*(g)))\in \texttt{\rm graph}(\partial_H^c \Psi).$ If
we set $\X\Psi(g)$ as $\X\psi_{g'}(g),$ for some $g'$ such that
$g\in \texttt{\rm dom}(\X\psi_{g'}),$ the optimal map $s_*$ is
given, almost surely, by
$$
s_*(g)=g\exp\left((\X c(g,\cdot))^{-1}(\X\Psi
(g))-\xi_1(g)\right).
$$
We conclude that the family of functions $\Psi$ plays the role of
the \lq\lq potential\rq\rq\ of the problem.

\bibliography{cacapi}

\bibliographystyle{plain}
\end{document}